\documentclass[journal]{IEEEtran}

\usepackage{amsmath,mathtools,amsfonts,amssymb,amsthm}
\usepackage{amsmath,amsfonts}
\usepackage{algpseudocode}
\usepackage{algorithm}
\usepackage{array}
\usepackage[caption=false,font=normalsize,labelfont=sf,textfont=sf]{subfig}
\usepackage{textcomp}
\usepackage{stfloats}
\usepackage{url}
\usepackage{verbatim}
\usepackage{graphicx}
\usepackage{soul}
\def\BibTeX{{\rm B\kern-.05em{\sc i\kern-.025em b}\kern-.08em
    T\kern-.1667em\lower.7ex\hbox{E}\kern-.125emX}}
\usepackage{balance}
\usepackage{multirow}
\usepackage{booktabs} 

\usepackage{caption}                         
\usepackage{algorithm}
\usepackage{algpseudocode}
\newtheoremstyle{sltheorem}
{}                
{}                
{}        
{10pt}                
{\bfseries}       
{:}               
{ }               
{}                
\theoremstyle{sltheorem}
\newtheorem{theorem}{Theorem}
\newtheorem{definition}{Definition}

\usepackage{xcolor}
\usepackage{cite}
\usepackage{enumitem,kantlipsum}

\begin{document}

\title{Efficient Scenario Generation for Chance-constrained Economic Dispatch Considering Ambient Wind Conditions}
\author{Qian Zhang,~\IEEEmembership{Student Member, IEEE,} Apurv Shukla,~\IEEEmembership{Member, IEEE,} and Le Xie,~\IEEEmembership{Fellow, IEEE}
\thanks{This work is supported in part by NSF ECCS-2038963, the U.S. Department of Energy (DoE) Office of Energy Efficiency and Renewable Energy (EERE) under the Solar Energy Technologies Office (SETO) Award Number DEEE0009031, and Texas A\&M Energy Institute. The authors are with the Department of Electrical and Computer Engineering, Texas A\&M University, College Station, TX, 77843 USA (e-mail: zhangqianleo@tamu.edu; apurv.shukla@tamu.edu; le.xie@tamu.edu).}}

\maketitle
\begin{abstract}
Scenario approach is an effective data-driven method for solving chance-constrained optimization while ensuring desired risk guarantees with a finite number of samples. Crucial challenges in deploying this technique in the real world arise due to {\color{black}non-stationarity environments and the absence of appropriate risk-tuning models tailored for the desired application.} In this paper, we focus on designing efficient scenario generation schemes for economic dispatch in power systems. We propose a novel scenario generation method based on filtering scenarios using ambient wind conditions. These filtered scenarios are deployed incrementally in order to meet desired risk levels while using minimum resources. In order to study the performance of the proposed scheme, we illustrate the procedure on case studies performed for both 24-bus and 118-bus systems with real-world wind power forecasting data. Numerical results suggest that the proposed filter-and-increment scenario generation model leads to a precise and efficient solution for the chance-constrained economic dispatch problem.
\end{abstract}

\section{Introduction}\label{sec:intro}
\indent 
{\color{black}The integration of intermittent renewable sources of generation into the existing power grids has posed significant challenge due to the inherent uncertainty associated with these sources. While both solar and wind generators introduce uncertainty, wind power, in particular, presents unique challenges. Unlike solar, wind energy lacks diurnal patterns, making it more difficult to predict. Furthermore, the scale of installed and consumed wind energy often surpasses that of solar power in many dispatch areas, leading to a higher level of uncertainty in the overall power system.} For instance, wind power accounts for 25\% of the total electricity consumed by the Electric Reliability Council of Texas (ERCOT) region during 2022, compared to $6\%$ for solar~\cite{ercot}.\\
\indent {\color{black}In response to these challenges, stochastic optimization (SO) and robust optimization (RO) are widely used methods for power engineers to deal with uncertainties \cite{roald2023power}. SO utilizes probabilistic models to handle randomness \cite{prekopa2013stochastic}, while RO employs set-based and deterministic uncertainty models \cite{ben2009robust}. We focus on chance-constrained optimization (CCO) which bridges the gap between probabilistic and deterministic approaches, providing explicit probabilistic guarantees on the feasibility of optimal solutions~\cite{geng2019data}.}\\ 
\indent Over the past decade, many attempts have been made at reformulating {\color{black}CCO} into a more {\color{black}computationally tractable} form. Bienstock \textit{et al.} \cite{bienstock2014chance} reformulate the chance-constrained DC-optimal power flow (DC-OPF) under affine control as a second-order cone program by a moment-based reformulation. Roald \textit{et al.} \cite{roald2017chance} extend similar reformulations into the AC-OPF model with some efficient algorithms. {\color{black} Other reformulations have also been proposed considering the distributionally robust property, especially in economic dispatch problem~\cite{zhang2016distributionally,poolla2020wasserstein}.}\\
\indent {\color{black}However, most of the above approaches require the explicit or approximated distribution function of uncertainty, which is hard to validate with streaming data. Data-driven optimization methods, unconstrained by specific underlying uncertainty distributions,} have received substantial attention in recent years~\cite{geng2019data}, especially the Sample Average Approximation (SAA)~\cite{luedtke2008sample} and the Scenario Approach~\cite{campi2008exact}. In the realm of power systems, applications of SAA are evident in day-ahead unit commitment~\cite{bagheri2017data}, capacity planning~\cite{madavan2022conditional}, and other offline domains~\cite{jirutitijaroen2008reliability}. In real-time cases, the scenario approach proves advantageous due to its rapid computational speed, 
avoiding the long solving time caused by binary variables in SAA. {\color{black}Despite many recent works that have tried to apply the scenario approach in power systems \cite{roald2023power,geng2019data,ming2017scenario,modarresi2018scenario,geng2021computing,vrakopoulou2013probabilistic,yan2022data,geng2019chance}, there still exist two main limitations.\\
\textit{A. Scenario Generation Limitation}\\
\indent Several studies validate scenario approach in the economic dispatch by using synthetic data created from representative distribution~\cite{ming2017scenario,modarresi2018scenario}. These techniques do not utilize and exploit the fact that the scenario approach can be distribution-agnostic since the scenarios can be directly extracted from previous experience. In most real-world settings, {\color{black}empirical} data is collected from non-stationary environments, wherein the distribution of the random variable depends on environmental conditions~\cite{campi2021scenario}. Directly sampling from the past will incorporate all environmental factors from {\color{black}empirical} data in scenario approach. For example, the wind power forecasting error addressed in this paper, is the primary source of uncertainty in the economic dispatch process ~\cite{bludszuweit2008statistical, hodge2012comparison}, while the wind forecast error distribution varies based on forecasting techniques, power output levels, and ambient conditions, as noted in studies by ~\cite{pinson2010conditional, zhang2013modeling}.\\ 
\textit{B. Risk Tuning Limitation}\\
\indent After choosing the scenario space, the next step is to decide the number of scenarios. The conventional \textit{sample and discard} approach \cite{campi2011sampling} requires the decision-makers to first create an estimate of the sample size and calculate the exact risk level ex-ante. Then, it is decided whether scenarios need to be discarded to trade off risk and performance. To ensure meeting the risk requirements, a conservative sample size is used, which could be astronomical for problems with a large number of decision variables. This inefficient risk-tuning process consumes superfluous data and results in a long computation time.\\
\indent The main contributions of this paper are dealing with these two limitations, which can be summarized as follows:\\
\begin{enumerate}[leftmargin=*]
    \item To enhance scenario accuracy, the conditional distribution is considered without assuming any knowledge of the true distribution.  Utilizing correlation analysis, wind forecast error scenarios are generated from empirical data similar to the present environment, ensuring a more precise representation of real-world conditions.\\
    \item The incremental risk tuning method is introduced to meet the risk requirement with minimum data resources~\cite{garatti2022complexity}. After declaring a desired risk level, scenarios are generated iteratively to eventually hit a desired level of risk.\\
    \item Algorithms are designed to incorporate scenario generation with risk-tuning processes efficiently.\\
\end{enumerate}}
\indent The remainder of this article is organized as follows: Section~\ref{sec:Problem} formulates the chance-constrained economic dispatch problem and highlights the challenges associated with solving the problem with conditional wind power forecast error. Section~\ref{sec:method} introduces our incremental scenario approach, and Section~\ref{sec:dis} discusses the assumptions and limitations of the proposed method. We demonstrate the efficacy of the proposed approach on 24 and 118-bus systems in Section~\ref{sec:case}.
\section{Problem Statements}\label{sec:Problem}
\subsection{Chance-Constrained Economic Dispatch}
We consider the chance-constrained DC-OPF formulation in the presence of wind-forecasting uncertainty~\cite{bienstock2014chance,vrakopoulou2013probabilistic}:
\begin{subequations} 
\label{cco}
\begin{align}
\min _{g, \eta}\; &c(g)   \label{subeqn:obj}\\
\text { s.t }&  \mathbf{1}^{\top} g=\mathbf{1}^{\top} d-\mathbf{1}^{\top} \hat{w} \label{subeqn:banl1}\\
& \underline{g} \preceq g \preceq \bar{g}\\
& f(\hat{w}, \tilde{w})=H_g\left(g-\mathbf{1}^{\top} \tilde{w} \eta\right)-H_d d+H_w(\hat{w}+\tilde{w}) \label{subeqn:f}\\
& \mathbb{P}_{\tilde{w}}\left(\begin{array}{l}\underline{f} \preceq f(\hat{w}, \tilde{w}) \preceq \bar{f}\\
\underline{g} \preceq g-\mathbf{1}^{\top} \tilde{w} \eta \preceq \bar{g}\\
R_d\preceq-\mathbf{1}^{\top} \tilde{w} \eta\preceq R_u \end{array}\right) \geq 1-\epsilon \label{subeqn:risk}\\
& \mathbf{1}^{\top} \eta=1 \label{subeqn:banl2}
\end{align}
\end{subequations}
The decision variables are generator output levels $g \in \mathbf{R}^{n_g}$, and an affine control policy $\eta \in \mathbf{R}^{n_g}$ {\color{black}proportionally allocating total wind fluctuation $\mathbf{1}^{\top} \tilde{w}$ to each generator}\footnote{\color{black} In this paper, $\mathbf{1}^{\top}$ denotes a row vector of all ones, with its dimension adjustable to match the vector it multiplies.}. The objective function is the total generations cost $c(g)$. The load level is $d \in \mathbf{R}^{n_d}$, and the wind generation $w=\hat{w}+\tilde{w}$ consists of \textit{deterministic} wind forecast value $\hat{w} \in \mathbf{R}^{n_w}$ and the \textit{uncertain} forecast error $\tilde{w} \in \Delta$, where $\Delta\subseteq\mathbf{R}^{n_w}$ is the uncertainty set. Transmission line flows $f \in \mathbf{R}^{n_f}$ are calculated using (\ref{subeqn:f}), where $H_g$, $H_d$, and $H_w $ are the corresponding sub-matrix of the power transfer distribution factor (PTDF) matrix $H$. Constraints include transmission line flow limits $[\underline{f}, \bar{f}]\in \mathbf{R}^{n_f} \times \mathbf{R}^{n_f}$, generation capacity limits $[\underline{g}, \bar{g}]\in \mathbf{R}^{n_g} \times \mathbf{R}^{n_g}$ and the ramp up(down) rate limits $[R_d, R_u]\in \mathbf{R}^{n_g} \times \mathbf{R}^{n_g}$ are modeled as a chance-constraint form under risk $\epsilon$ in (\ref{subeqn:risk}). \\ 
\indent As mentioned in~\cite{vrakopoulou2013probabilistic}, the affine control policy $\eta$ only focuses on the steady-state behavior of the \textit{Automatic Generation Control} (AGC) action in dispatch time scale, i.e. $5$ to $15$ minutes, but not the $2$ to $6$ seconds fast time-scale regulation process. The system's active power deviation is allocated to generators based on $\eta$, {\color{black}which is also well known as \textit{participation factors} in the conventional AGC scheme \cite{kundur2022power}. It is easy to confirm that constraints (\ref{subeqn:banl1}) and (\ref{subeqn:banl2}) imply the supply-demand balance in the presence of wind uncertainties:
\begin{equation} \label{balance}
\mathbf{1}^{\top} (g -\mathbf{1}^{\top} \tilde{w} \eta)=\mathbf{1}^{\top} d-\mathbf{1}^{\top} (\hat{w}+\tilde{w})
\end{equation}}
\indent Setting the constant {\color{black}affine control policy} $\eta$ prior to the next dispatch interval will unavoidably be economically inefficient if the net load's fluctuation or forecasting error is large. Changing the $\eta$ more frequently within dispatch interval, or incorporating the optimization program into the AGC control policy may improve the economic efficiency \cite{liu2012enhanced,li2015connecting}, but these methods ignore the network constraints and are hard to apply to the bulk power system due to the communication delay or solving time. In this paper, we mainly concentrate on improving the dispatch performance in the chance-constrained problem. The description of {\color{black}affine control policy} $\eta$ in (\ref{cco}) is based on two assumptions: 1) All the traditional generators participate in the AGC actions; 2) The whole system is regarded as one control area.

\subsection{Conditional Wind Power Forecast Error}
Due to the nonlinear wind power curve, the wind power forecast error is observed to vary with the level of its output~\cite{lange2005uncertainty,lei2009review}, while spatial and weather parameters also indirectly affect the forecasting quality ~\cite{miettinen2020simulating,hanifi2020critical}. Based on these facts, the wind power forecast error measure $\mathbb{P}_{\tilde{w}}$ in (\ref{subeqn:risk}) should be modeled as a conditional probability distribution from similar environments. Previous studies have focused on approximating the measure of conditional forecast error, eg., ~\cite{lange2005uncertainty} generate the conditional error model based on the wind turbine power curve, ~\cite{pinson2010conditional} employ a fuzzy inference model to obtain conditional prediction intervals and~\cite{zhang2013modeling} calculate the conditional forecasting error from joint distributed data using copula theory. All the above approximation methods for the distribution $\mathbb{P}_{\tilde{w}}$ under similar conditions are incompatible with the data-driven approach in chance-constraint optimization, which uses distribution agnostic {\color{black}empirical} data. Furthermore, the wind power output, temperature, and weather parameters are continuous variables, meaning the past scenarios' observations will be distinct with probability 1, which makes it impossible to generate scenarios in the data-driven scheme from an identical environment but similar environments. Scenario selection from similar environments presents a major engineering obstacle: a large search space results in conservative decision-making while a smaller search space results in limited data~\cite{campi2021scenario}.

\section{Method}\label{sec:method}
\subsection{Scenario Approach}
\indent The scenario approach randomly extracts $N$ independent and identically distributed (\textit{i.i.d.}) scenarios to approximate the chance-constrained program. Supposing we have the random wind forecasting error scenarios set $\mathcal{N}:=\left\{\tilde{w}_1, \tilde{w}_2, \cdots, \tilde{w}_N\right\}$, the chance-constrained inequalities (\ref{subeqn:risk}) in DC-OPF problem can be replaced by scenario-based inequalities (\ref{sa}):
\begin{subequations} \label{saaa}
\begin{align}
\min _{g, \eta}\; &c(g)   \label{subeqn:obj}\\
\text { s.t }&  \mathbf{1}^{\top} g=\mathbf{1}^{\top} d-\mathbf{1}^{\top} \hat{w} \label{subeqn:banl1}\\
& \underline{g} \preceq g \preceq \bar{g}\\
& f(\hat{w}, \tilde{w})=H_g\left(g-\mathbf{1}^{\top} \tilde{w} \eta\right)-H_d d+H_w(\hat{w}+\tilde{w}) \label{subeqn:f}\\
&\begin{array}{l}\underline{f} \preceq f(\hat{w}, \tilde{w}_i) \preceq \bar{f}\\
\underline{g} \preceq g-\mathbf{1}^{\top} \tilde{w}_i \eta \preceq \bar{g}\\
R_d\preceq-\mathbf{1}^{\top} \tilde{w}_i \eta\preceq R_u \end{array} \quad i=1,2,3,...,N \label{sa}\\ 
& \mathbf{1}^{\top} \eta=1 \label{subeqn:banl2}
\end{align}
\end{subequations}

\indent To distinguish from the original optimization problem (\ref{cco}), we name the above \textit{scenario problem} as $\mathrm{SP}(\mathcal{N})$.
\begin{definition}[Violation Probability]
The \textit{vilolation probability} of a candidate solution $(g^*,\eta^*)$ is defined as the probability that $(g^*,\eta^*)$ is infeasible, i.e., $\mathbb{V}_{\tilde{w}}(g^*,\eta^*):=\mathbb{P}_{\tilde{w}}((g^*,\eta^*)\notin \mathcal{X}_{\tilde{w}})$, where $\mathcal{X}_{\tilde{w}}$ is the decision set generated by $\mathrm{SP}(\mathcal{N})$.
\end{definition}
\begin{definition}[Support Constraint]
The scenario-dependent constraint corresponding to sample $\tilde{w}_s, s \in \{1,2,...,\mathcal{S}\}$, is a \textit{support constraint} or \textit{support scenario} if its removal improves the solution of $\mathrm{SP}(\mathcal{N})$, i.e., if it decreases the optimal cost (\ref{subeqn:obj}).    
\end{definition}
\begin{definition}[Sample Complexity]
The number of support scenarios in $\mathrm{SP}(\mathcal{N})$ is defined as the sample complexity.    
\end{definition}
\begin{definition}[Helly's Dimension]
Helly's dimension of the scenario problem $\mathrm{SP}(\mathcal{N})$ is the smallest integer $h$ that $h\geq_{\text{ess sup}\mathcal{N}\subseteq\Delta^N}|\mathcal{S(N)}|$ holds for any finite $N\geq1$, where $|\mathcal{S(N)}|$ is the number of support constraints or sample complexity \footnote{{\color{black}Because the randomness of sampling, the number of support constraints might be different especially when sample size is small. Here $_\text{ess sup}$ means the essential supremum to ignore some exceptional cases.}}.    
\end{definition}
\indent The most important contribution of the scenario approach is the relationship between violation probability $\mathbb{V}_{\tilde{w}}(g^*,\eta^*)$, the number of scenarios $N$ and the sample complexity.
\begin{theorem}[Exact Feasibility \cite{campi2008exact,campi2018general}] \label{theorem1}
Under the assumptions of nondegeneracy and feasibility of the optimization problem, the deepest results show that the distribution of $\mathbb{V}_{\tilde{w}}(g^*,\eta^*)$ is dominated by a Beta distribution, namely:
\begin{equation}
\label{ExactSP}
\mathbb{P}_{\tilde{w}}^N\left(\mathbb{V}_{\tilde{w}}(g^*,\eta^*)>\epsilon\right) \leq \sum_{i=0}^{h-1}\left(\begin{array}{c}
N \\
i
\end{array}\right) \epsilon^i(1-\epsilon)^{N-i} := \beta
\end{equation}
\end{theorem}
where $h$ is the Helly's dimension of $\mathrm{SP}(\mathcal{N})$, and $1-\beta$ is defined as the confidence bound for the solution based on any $N$ \textit{i.i.d} scenarios.
\begin{theorem}[Property for Convex Problem \cite{campi2008exact}] \label{hn}
Supposing all the constraints in $\mathrm{SP}(\mathcal{N})$ is convex for every instance of $\tilde{w}$, the sample complexity $|\mathcal{S(N)}|$ for $\mathrm{SP}(\mathcal{N})$ is at most $n$. In other words, $h \leq n$, where $n$ is the number of decision variables after eliminating the equality constraints {\color{black}\cite{boyd2004convex}}.\\
\end{theorem}
\indent For the convex problem, Helly's dimension $h$ can be replaced by $n$ to simplify the problem by applying Theorem \ref{hn}, but it often causes extremely conservative results. To compute the lower bound of $h$, we suggest using the dual-based Algorithm~\ref{alg:dual} proposed in~\cite{geng2021computing}.
\begin{algorithm}[H]
\caption{Searching Support Scenarios Using Dual Variables}\label{alg:dual}
Solving the \textit{scenario problem} $\mathrm{SP}(\mathcal{N})$ \\
Generate the primal solution $(g^*_\mathcal{N},\eta^*_\mathcal{N})$ and the constraints (\ref{sa}) related dual solution $\mu^*_i$ ($i=1,2,3,...,N$) \\
Let $\mathcal{M} = {i \in \mathcal{N}:\|\mu^*_i\|>0}$. Set $\mathcal{S} \gets \O$
\begin{algorithmic}
\For{$i \in \mathcal{M}$} 
\State Solve $\mathrm{SP}_{\mathcal{M}-i}$ and compute $(g^*_{\mathcal{M}-i},\eta^*_{\mathcal{M}-i})$
\If{$c(g^*_{\mathcal{M}-i})<c(g^*_\mathcal{N})$}
\State $\mathcal{S} \gets \mathcal{S} + i$
\EndIf
\EndFor
\end{algorithmic}
\textbf{OUTPUT:} The support scenarios $\mathcal{S}$
\end{algorithm}

\subsection{Sampling Scenarios from Parameter Space}
\indent Many environmental parameters, such as location, wind speed, temperature, wind direction, and relative humidity affect the wind power forecasting quality \cite{hanifi2020critical}. In this paper, we suppose the wind generators are from the same area with the same forecasting algorithm, while the \textit{deterministic} wind forecast value, the wind power changing rate, temperature, and relative humidity are selected as the four main parameters that affect the forecasting error.\\
\indent \textit{Remark}: The \textit{deterministic} wind forecast value and the wind power ramp rate are parameters that are integrated with other environmental information, especially the wind speed and wind speed ramp rate, which are observed having a close relationship with the wind power forecasting quality.
\begin{definition}[Parameter Space for Wind Power Forecasting Error]
The parameter space $\mathcal{V}_\mathcal{N}$ is defined as the set of environmental parameters which the past scenarios $\mathcal{N}$ are extracted from. For instance, the temperature between $70^\circ F$ and $80^\circ F$ is a temperature parameter space. 
\end{definition}
\begin{definition}[Probability Distribution Over Parameter Space]
Let $\mathbb{P}_{\tilde{w}|\mathcal{V}_\mathcal{N}}$ be a probability distribution over the parameter space $\mathcal{V}_\mathcal{N}$. 
\end{definition}
\indent \textit{Remark}: Strictly speaking, $\mathcal{V}_\mathcal{N}$ should be identical with the environment parameter at the forecasting moment $\mathcal{V}_{\text{Now}}$ to acquire more precise risk guarantee. However, the input data in the scenario approach is directly extracted from the empirical experience with continuous environment parameters. Finding the past scenarios under the same environment parameter as the future is impossible with probability 1. Even in some frontier probabilistic prediction methods, it is also hard to guarantee the accuracy of predicted conditional distribution but check how close it approximates to the real distribution based on testing data \cite{vovk2012conditional}.\\
\indent Based on the definitions above, the barrier of bringing the scenario approach to the real world is to find the parameter space $\mathcal{V}_\mathcal{N}$ which includes both the potential environment parameter in the future and enough number of empirical scenarios to meet the risk requirement in Theorem \ref{theorem1}.
\subsection{Main Result I: Correlation-Based Scenario Generation}
\begin{figure*}
\centering
  \includegraphics[width=0.8\textwidth]{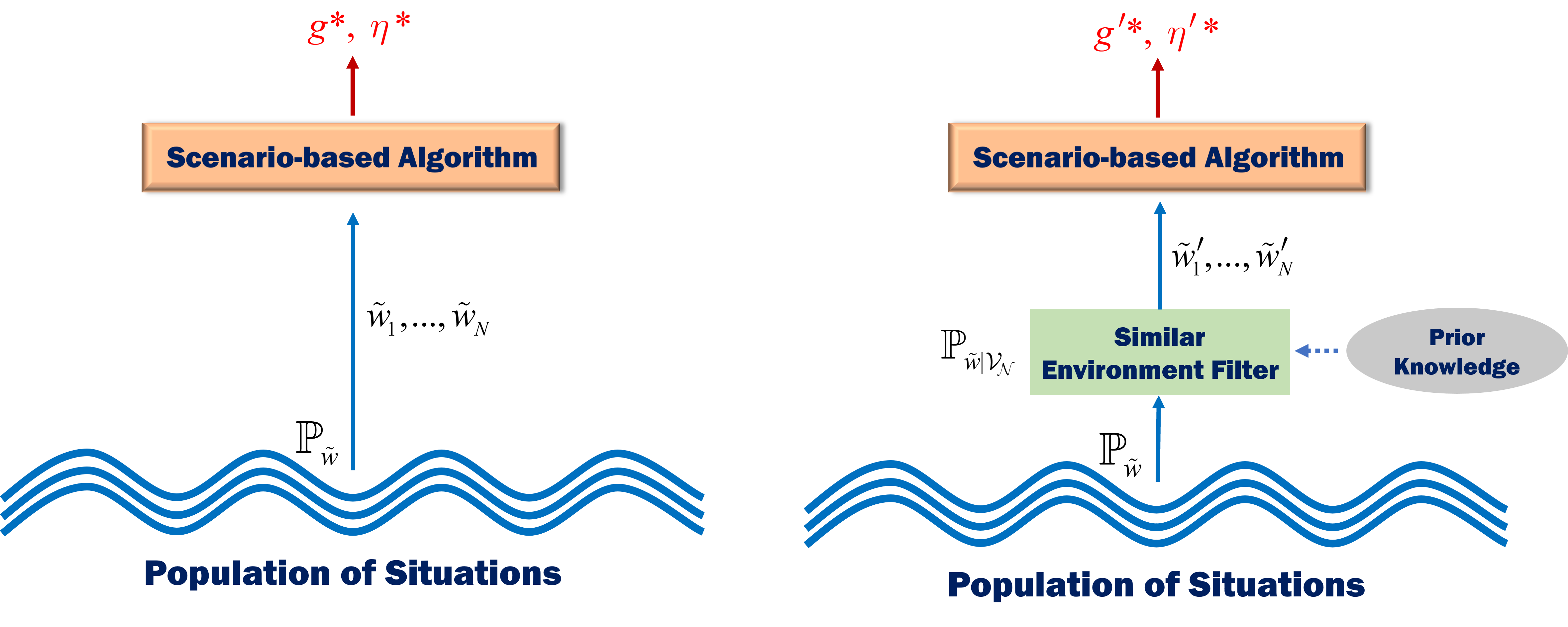} 
  \caption{The comparison of the conventional scenario approach (left) and the proposed scenario generation model (right)}
  \label{scheme}
\end{figure*}

\indent The Pearson correlation coefficient is used to quantify the relationship between each environmental parameter and the wind power forecasting error from the past forecasting data, for example, the past half-year data before the decision-making day. Let random variable $P$ represent one of the environmental parameters, and random variable $\tilde{W}$ denote the wind power forecasting error, then the Pearson correlation coefficient between these two variables is given by:
\begin{equation}
\label{correlation1}
\rho_{P, \tilde{W}}=\frac{\mathbb{E}\left[(P-\mu_P)(\tilde{W}-\mu_{\tilde{W}})\right]}{\sigma_P \sigma_{\tilde{W}}}
\end{equation}
where $\sigma_P$ and $\sigma_{\tilde{W}}$ are the standard deviation of $P$ and $\tilde{W}$, while $\mu_P$ and $\mu_{\tilde{W}}$ are the mean value of $P$ and $\tilde{W}$, respectively. After calculating the correlation coefficient of each parameter, we combine the penalized parameters as the \textit{indicator vector} for estimating the difference in environmental {\color{black}conditions}. Suppose $\rho_1,\rho_2,\rho_3,\rho_4$ are the correlation coefficients of the four normalized parameters $p_1, p_2, p_3, p_4$ that affect the past forecasting error respectively, then the \textit{indicator vector} $v$ is constructed as:
\begin{equation}
\label{correlation2}
v=[\rho_1p_1, \rho_2p_2, \rho_3p_3, \rho_4p_4]^T
\end{equation}
\indent Based on the \textit{indicator vector}, we can define the wind power forecasting environment difference $d_{ij}$ as the distance between two forecasting environments $i$ and $j$, that is:
\begin{equation}
\label{correlation3}
d_{ij}=\|v_i-v_j\|
\end{equation}
{\color{black} where $\|\cdot\|$ can be any norm distance and Euclidean norm is used in the case study part.}\\
\indent After calculating the environment difference $d_{ij}$ between the now and the past, we can pick scenarios from the empirical data under more similar decision-making environments. The overview of {\color{black}the} proposed scheme and the conventional scenario approach is drawn in Fig.\ref{scheme}.\\
\indent Except for correlation analysis, many statistical or learning methods can be used as the similar environment filter represented in Fig.\ref{scheme}, such as the coefficient of determination ($R^2$), and other learning-based ranking algorithms \cite{cao2007learning}. Furthermore, some generative models, such as generative adversarial networks (GAN), are another path to get the scenarios under a similar environment. In this paper, we do not compare the accuracy of different filters but rather focus on illustrating this scheme and its performance on the economic dispatch problem. 
\subsection{Main Result II: Incremental Scenario Optimization}
\indent To find the proper number of scenarios meeting the risk requirement, two main scenario-based algorithms exist for convex problems. The simple one is called \textit{a-priori} approach \cite{campi2008exact}, where Helly's dimension $h$ is approximated by its upper bound: the number of decision variables $n$. In a-priori approach, the number of the needed scenarios can be directly given by Theorem \ref{theorem1} before solving the optimization program. The scenario sampling from a similar environment embedded in a-priori scenario optimization approach is summarized in the Algorithm \ref{alg:samp}.\\
\indent The a-priori approach has less computational complexity and risk guarantee before solving the optimization problem. When the number of decision variables is small, the a-priori approach is an efficient choice, but as the size of the system increases it leads to extreme conservatism. For instance, the sample complexity of the look-ahead economic dispatch problem~\cite{modarresi2018scenario} is $3-5$ while the decision variables $n$ is 864.\\
\begin{algorithm}[H]
\caption{A-priori scenarios approach sampling from similar environment}\label{alg:samp}
\textbf{INPUT1:} The past forecasting error data with environment parameters\\ 
\textbf{INPUT2:} The environment parameter at decision-making time 
\textbf{INPUT3:} The risk and confidence level $\epsilon$ and $\beta$
\begin{algorithmic}
\State $N \gets (\ref{ExactSP})$
\If{$N >$ empirical data size}
    \State Reset risk and confidence level
\Else {\quad sampling the scenarios $\mathcal{N}$ under similar environments}
    \State $d_{ij} \gets (\ref{correlation1})(\ref{correlation2})(\ref{correlation3})$
    \State $\mathcal{N} \gets$ the scenarios with the $N$ smallest $d_{ij}$
\EndIf
\State Solve the scenario-based optimization problem $\mathrm{SP}(\mathcal{N})$ with optimal solution $(g^*,\eta^*)$.
\end{algorithmic}
\textbf{OUTPUT:} The optimal solution $(g^*,\eta^*)$
\end{algorithm}
\indent Leveraging on support scenario searching algorithms, such as Algorithm \ref{alg:dual}, the a-priori approach can be improved if the risk level is updated after solving the optimization problem, called the \textit{a-posteriori} approach \cite{campi2021scenario}. In many practical engineering situations, the decision makers are interested in seeking the optimal solution given a certain risk level, implying that if the updated risk level is too conservative, the decision maker trades off risks for better performance by discarding some scenarios~\cite{campi2011sampling}. The computation complexity of this risk-tuning process can be reduced by using an incremental scenario optimization algorithm~\cite{garatti2022complexity}. In this paper, we propose an improved incremental scenario optimization algorithm based on sampling from a similar environment, see Algorithm~\ref{alg:incre}. In practice, this algorithm can generate the optimal solution with the given risk guarantee and provide some higher-risk solutions with better economic performance to meet the risk-tuning {\color{black}needs} of system operators.
\begin{algorithm}[H]
\caption{Incremental scenario optimization sampling from similar environment}\label{alg:incre}
\textbf{INPUT1:} The past forecasting error data with environment parameters\\ 
\textbf{INPUT2:} The environment parameter at decision-making time 
\textbf{INPUT3:} The risk and confidence level $\epsilon$ and $\beta$
\begin{algorithmic}[1]
\State Set $j:=1$ and $N_0=0$
\State Suppose the problem has $j$ support scenarios, i.e. $h=j$, and calculate the number of needed scenarios $N_j$ from (\ref{ExactSP}).
\State Collect a sample of scenarios $\tilde{w}_{N_{j-1}+1},\tilde{w}_{N_{j-1}+2},...,\tilde{w}_{N_{j}}$ with $N_{j}-N_{j-1}$ smallest environment difference $d$ without replacement.
\State Solve the scenario-based optimization problem $\mathrm{SP}(\mathcal{N}_j)$ with optimal solution $(g^*,\eta^*)_{N_j}$.
\State Compute the exact sample complexity $h_j$ of the solution $(g^*,\eta^*)_{N_j}$ based on Algorithm \ref{alg:dual}.
\State (Optional) Compute the exact risk level $\epsilon_j$ corresponding to $(g^*,\eta^*)_{N_j}$ after updating $h=h_j$ in (\ref{ExactSP}).
\If{$h_j \leq j$}
    \State halt the algorithm and \Return $(g^*,\eta^*):=(g^*,\eta^*)_{N_j}$
    \Else{\quad set $j:=j+1$} and GOTO step 2
\EndIf
\end{algorithmic}
\textbf{OUTPUT1:} The optimal solution $(g^*,\eta^*)$.\\
\textbf{OUTPUT2:} (Optional) The solution $(g^*,\eta^*)_{N_j}$ with higher risk level $\epsilon_j$.
\end{algorithm}
\indent In the previous sample and discard risk tuning method \cite{campi2011sampling}, the sample complexity is first supposed to be its upper bound, i.e. the number of decision variables, which means a large sample size may be used to solve the problem. After calculating the true sample complexity and risk level based on the solution, the {\color{black}decision-makers} can trade off the risk and performance by gradually discarding some support scenarios. Compared {\color{black}to} the incremental optimization method, the sample and discard method tunes the risk in a decremented way, which is inefficient when the true sample complexity is much smaller than the number of decision variables.

\section{Discussion}\label{sec:dis}
\indent The most important assumption in the scenario approach is the samples are identical independent distributed (\textit{i.i.d.}) random variables. In this section, we will discuss how the proposed method meets this \textit{i.i.d.} property and our weakness.   
\subsection{Identical Distribution}
\indent In the chance-constrained problem, the uncertainties are modeled as random variables sampled from the identical probability distribution. For example, the uncertainty of wind forecasting error in chance-constrained DC-OPF problem (\ref{cco}) is sampled from the fixed probability distribution $\mathbb{P}_{\tilde{w}}$. Because each forecasting happens in a complex and varying environment, the precise description of uncertainty should be under the environment condition, i.e. $\mathbb{P}_{\tilde{w}|\mathcal{V}_{\text{Now}}}$. For data-driven optimization methods, if we apply the whole empirical data to describe uncertainty, the identical probability distribution $\mathbb{P}_{\tilde{w}}$ should be regarded as the marginal distribution over the (past) environment parameters.\\
\indent Existing papers focusing on modeling the gap between the solution under empirical distribution $\mathbb{P}_{\tilde{w}}$ and exact distribution $\mathbb{P}_{\tilde{w}|\mathcal{V}_{\text{Now}}}$ require prior knowledge of the distance or mean value between two probability measures, which is hard to calculate in the wind forecasting scenario. Some machine learning approaches \cite{wan2013probabilistic} may be useful to approximate this exact distribution, but their results are not compatible with the direct data-driven program, especially the scenario approach.\\
\indent In this paper, instead of quantifying the gap between empirical distribution and exact distribution, 
we seek to filter empirical data through the similar environment parameter space $\mathcal{V}_\mathcal{N}$. The proposed solution's risk guarantee is based on the conditional distribution $\mathbb{P}_{\tilde{w}|\mathcal{V}_\mathcal{N}}$, but the simulations in the next section show that the testing results also meet the setting risk threshold even in the real-world data.
\subsection{Independent Random Variables}
\indent The purpose of sampling scenarios from parameter space is to make the conditional distribution $\mathbb{P}_{\tilde{w}|\mathcal{V}_\mathcal{N}}$ much closer to the exact one $\mathbb{P}_{\tilde{w}|\mathcal{V}_{\text{Now}}}$, which do not affect the independent property of scenarios under the two assumptions. The first is that the forecasting algorithm does not use the previous forecasting error data, which is true for most updated forecasting methods \cite{hanifi2020critical}. The second assumption is that the dispatch decision itself will not affect {\color{black}future forecasting errors}, which is also justifiable because the dispatch decision does not interact with the whole weather system in the short term. 
\section{Case Study}\label{sec:case}
\begin{figure}[H]
\centering
  \includegraphics[scale=0.6]{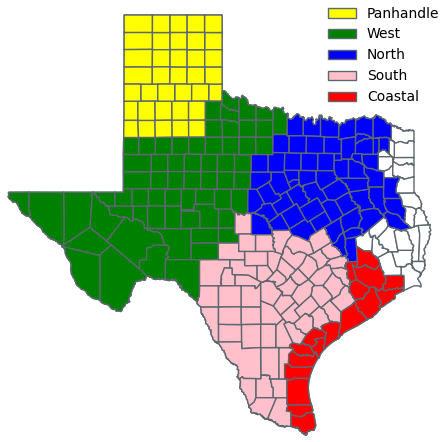}
  \caption{{\color{black}The wind forecasting region in ERCOT}}
  \label{forecastregion}
\end{figure}

\indent The scenario approach formulation (\ref{ExactSP}) has been validated in the chance-constrained economic dispatch many times \cite{geng2019data,ming2017scenario,modarresi2018scenario}. The scenarios in these previous papers were all obtained by sampling from some particular distribution, such as normal and beta distribution, which didn't utilize and exploit the distribution-free advantage of the original scenario theory.\\
\indent In our simulation, the knowledge about uncertainty is acquired directly through experience, i.e. the past recorded data. We focus on the economic dispatch under the 5-min unit, where the 5-min ahead wind forecast uncertainty plays an important part. The peak hours 16:00-18:00 in August 2022 is selected as the testing {\color{black} period} with 744 dispatch intervals. For each testing interval, the empirical forecasting error scenarios are generated from the past three or six months, i.e. a dataset with 25920 or 51840 5-min ahead forecasting errors with their corresponding environment parameters. All the data is acquired from the ERCOT website {\color{black} based on the five wind forecasting regions in Fig.\ref{forecastregion} } \cite{ercot}, while the weather data is from \cite{weather}.
\indent {\color{black}All the problems are solved using 64 GB RAM on the Intel XEON-10885M CPU (2.4GHz). The mathematical models were formulated using YALMIP on Matlab R2023a and solved using Gurobi v9.5.} 
\subsection{Conditional Wind Forecasting Error}
\begin{figure}[H]
\centering
  \includegraphics[scale=0.3]{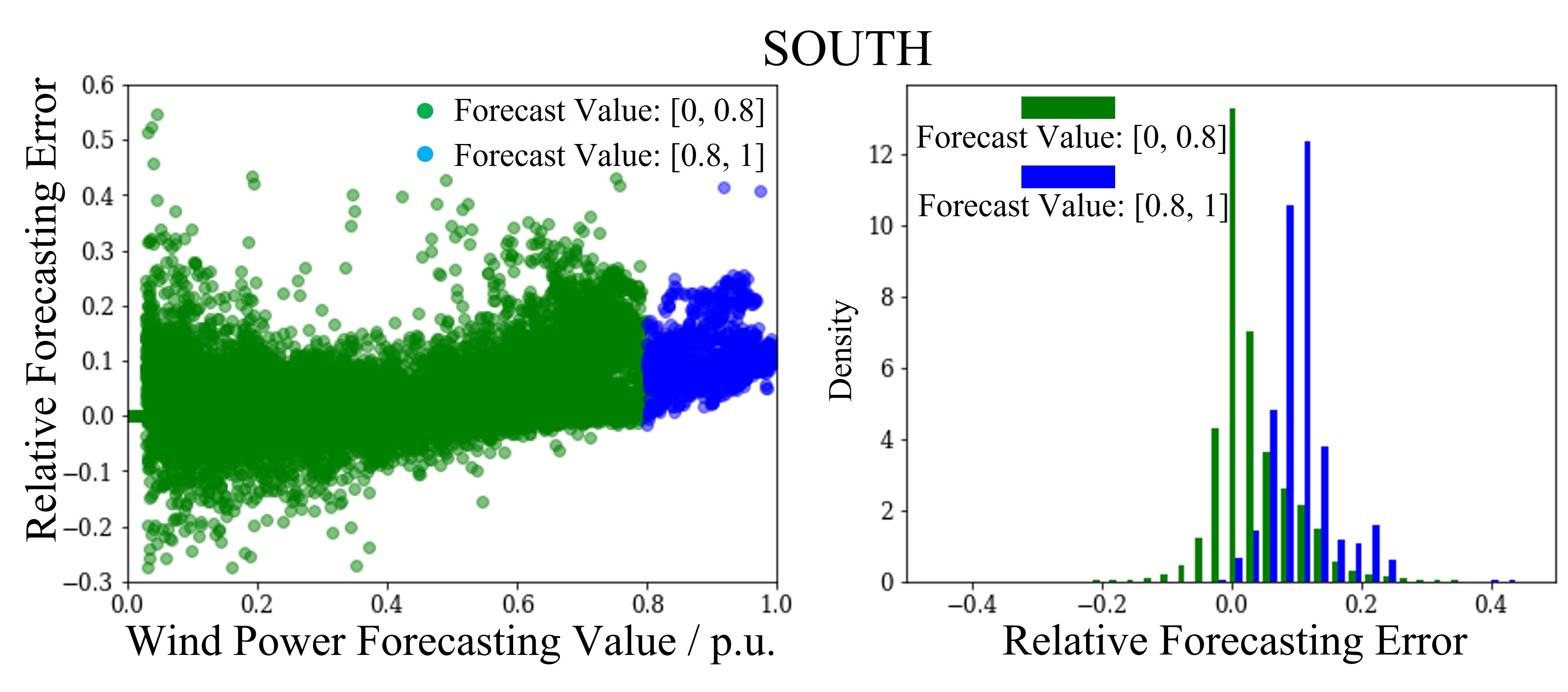}
  \includegraphics[scale=0.3]{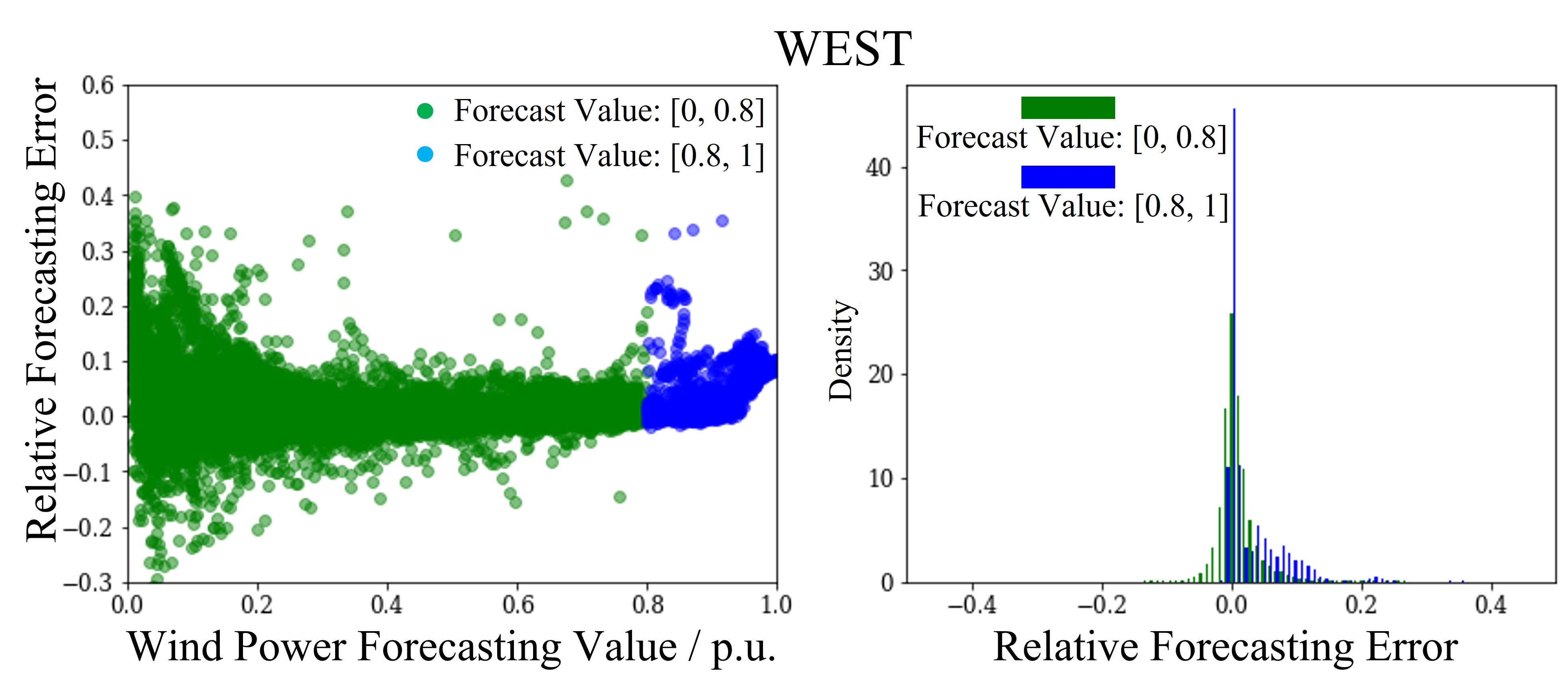}
  \caption{{\color{black}The scatter plot (left) and density function (right) of forecasting error under different wind power forecasting value}}
  \label{density_power}
\end{figure}

\begin{figure}[H]
\centering
  \includegraphics[scale=0.30]{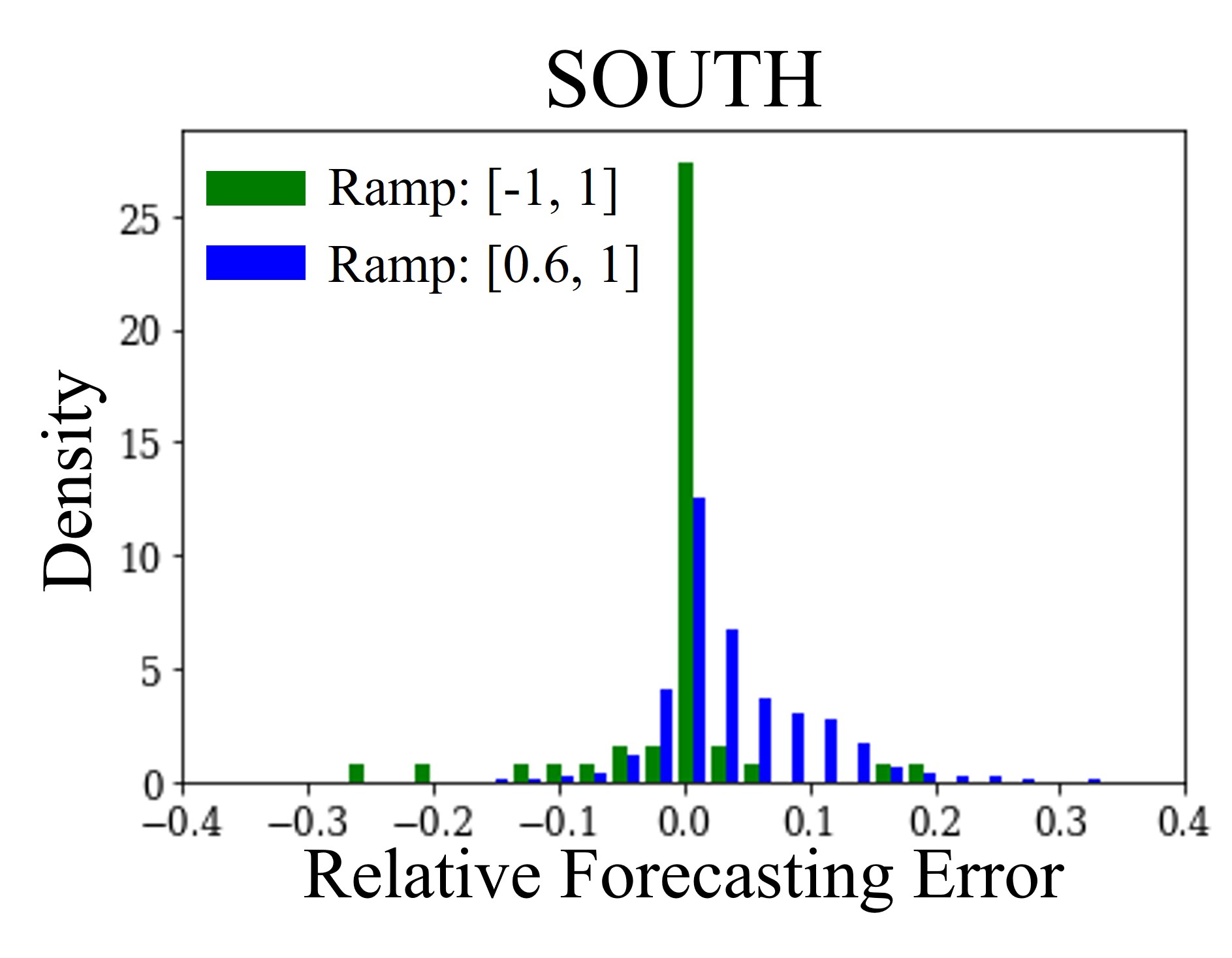}
  \includegraphics[scale=0.30]{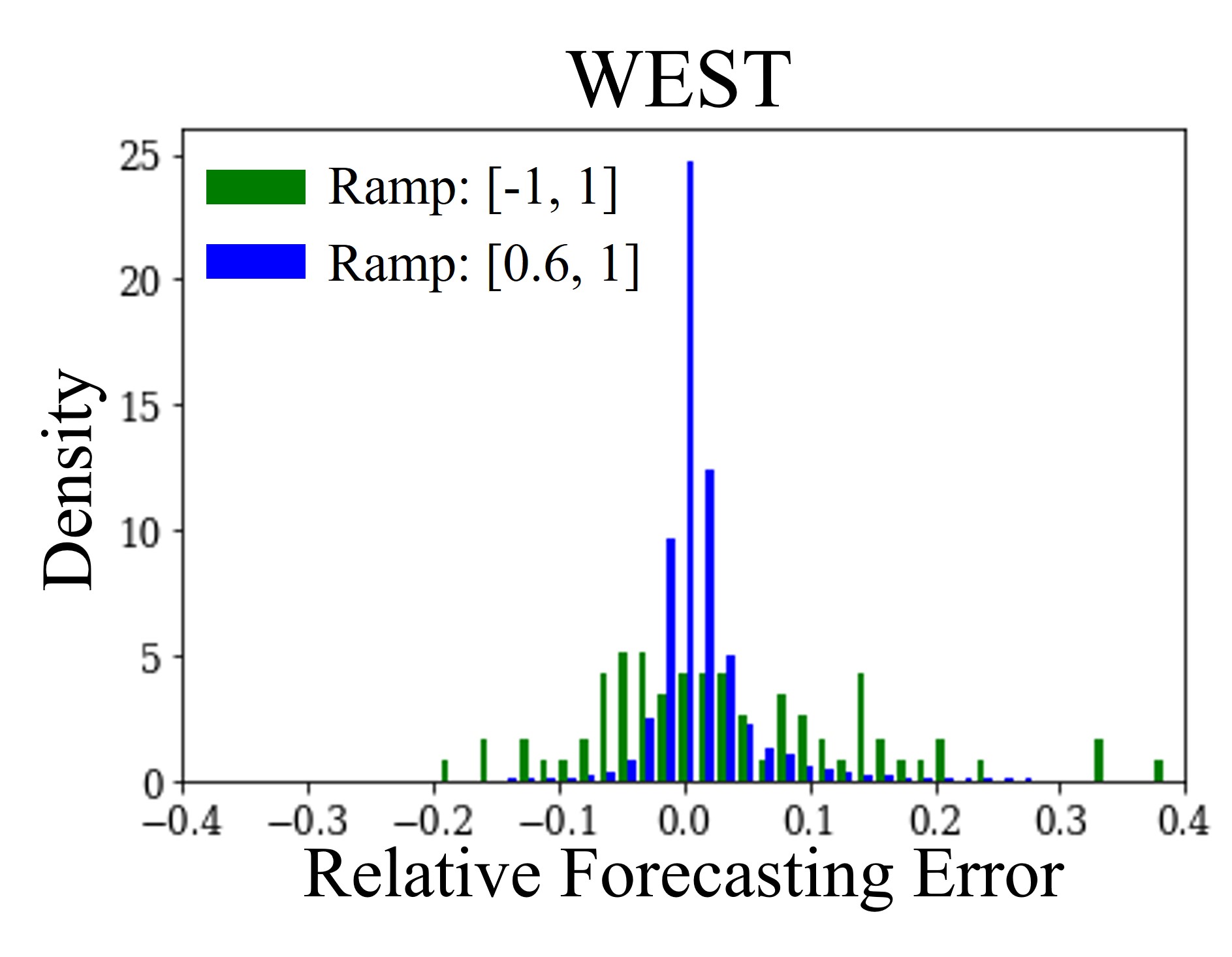}
  \caption{{\color{black}The density function of forecasting error under different wind power ramping rate}}
  \label{density_rate}
\end{figure}

\begin{figure}[H]
\centering
  \includegraphics[scale=0.30]{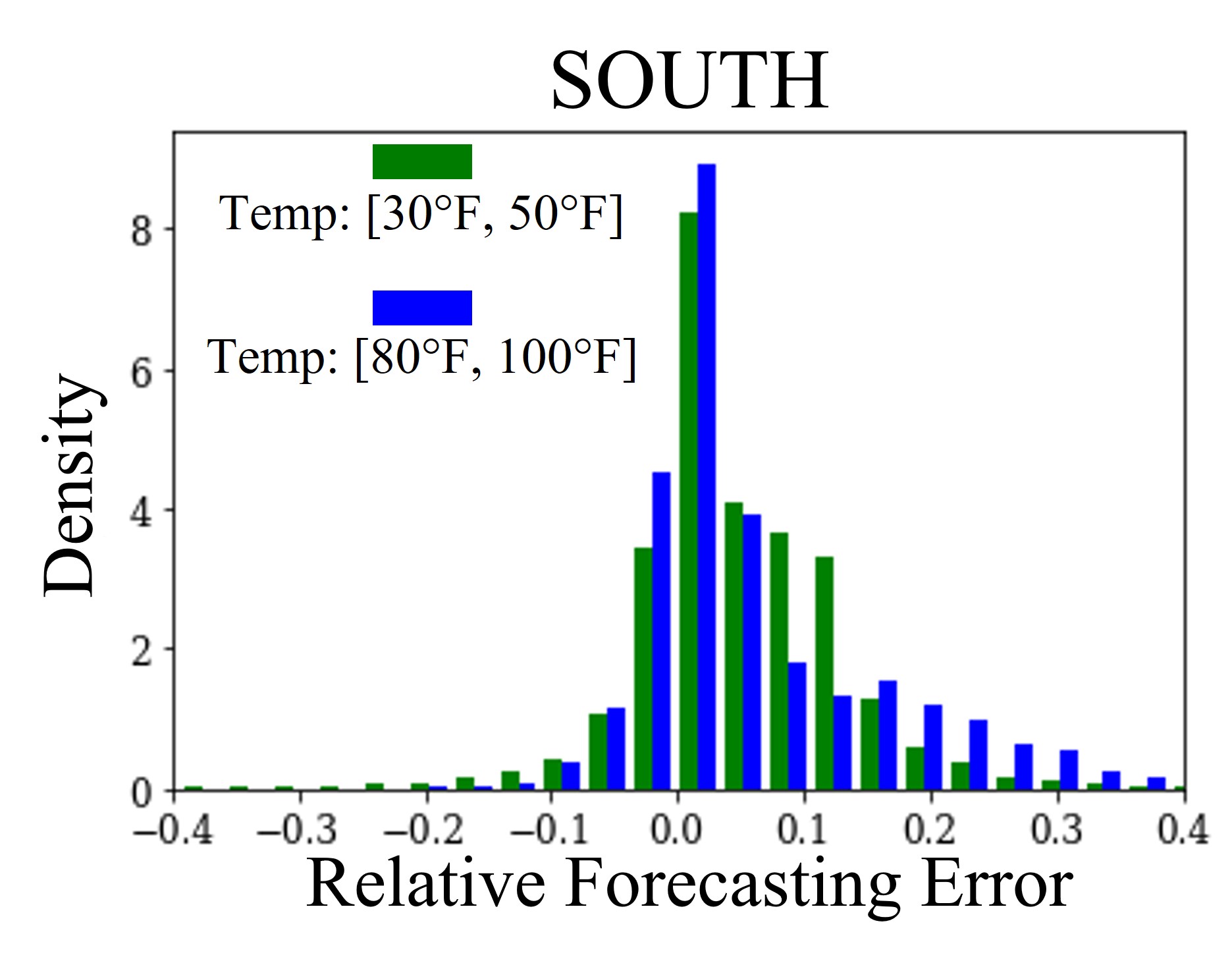}
  \includegraphics[scale=0.30]{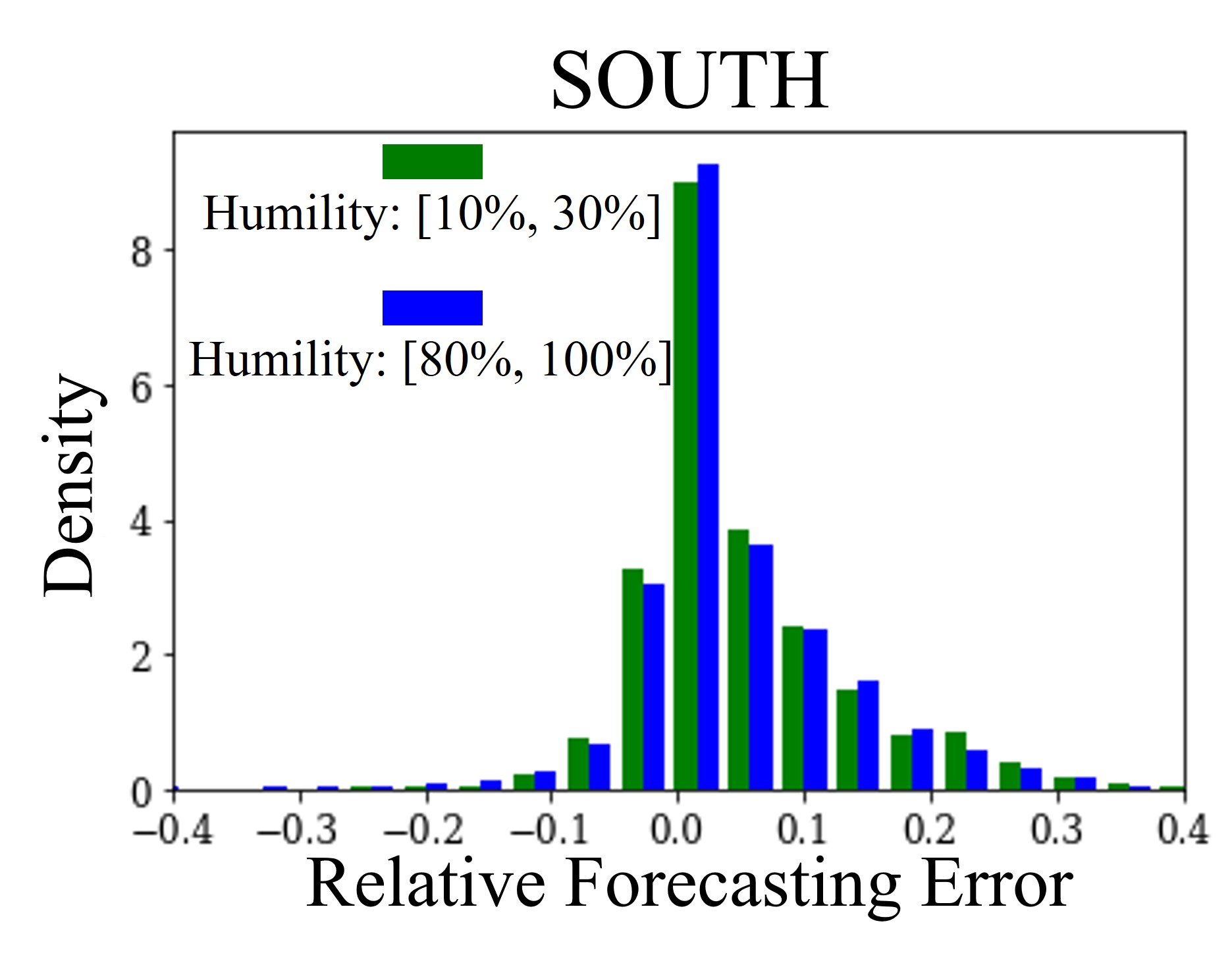}
  \caption{{\color{black}The density function of forecasting error under different temperatures and relative humidity}}
  \label{density_others}
\end{figure}

\begin{figure}[H]
\centering
  \includegraphics[scale=0.24]{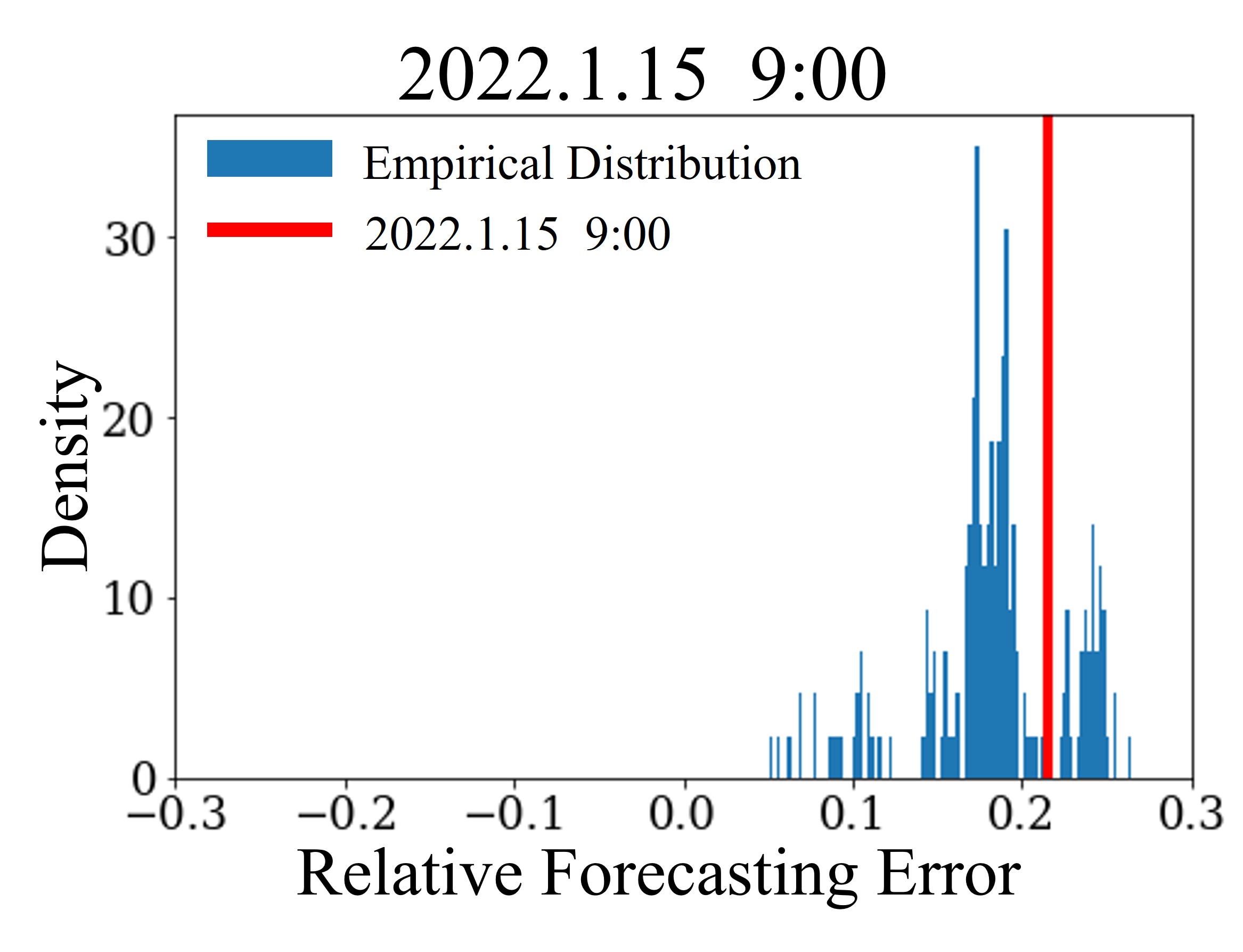}
  \includegraphics[scale=0.24]{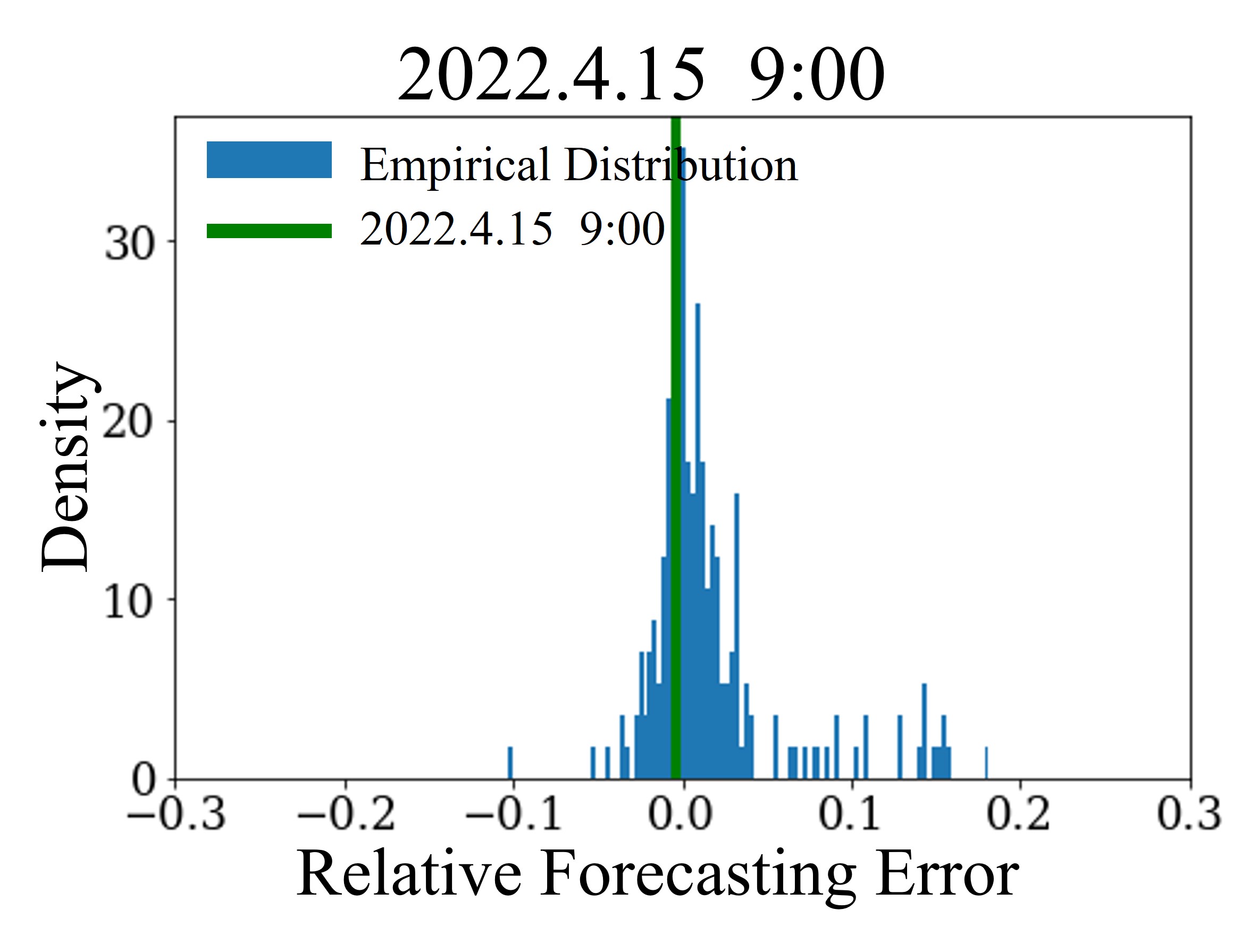}
  \includegraphics[scale=0.24]{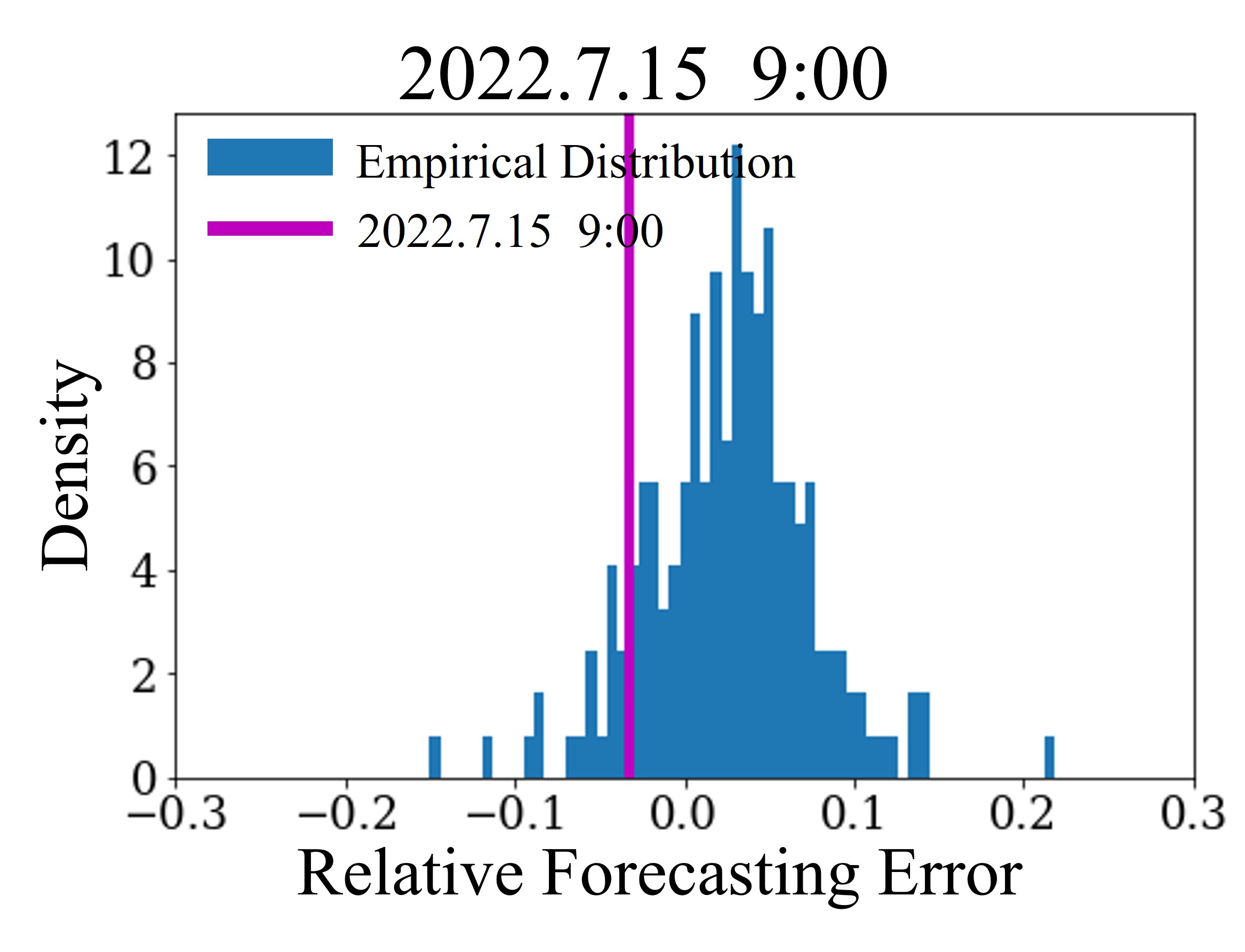}
  \includegraphics[scale=0.24]{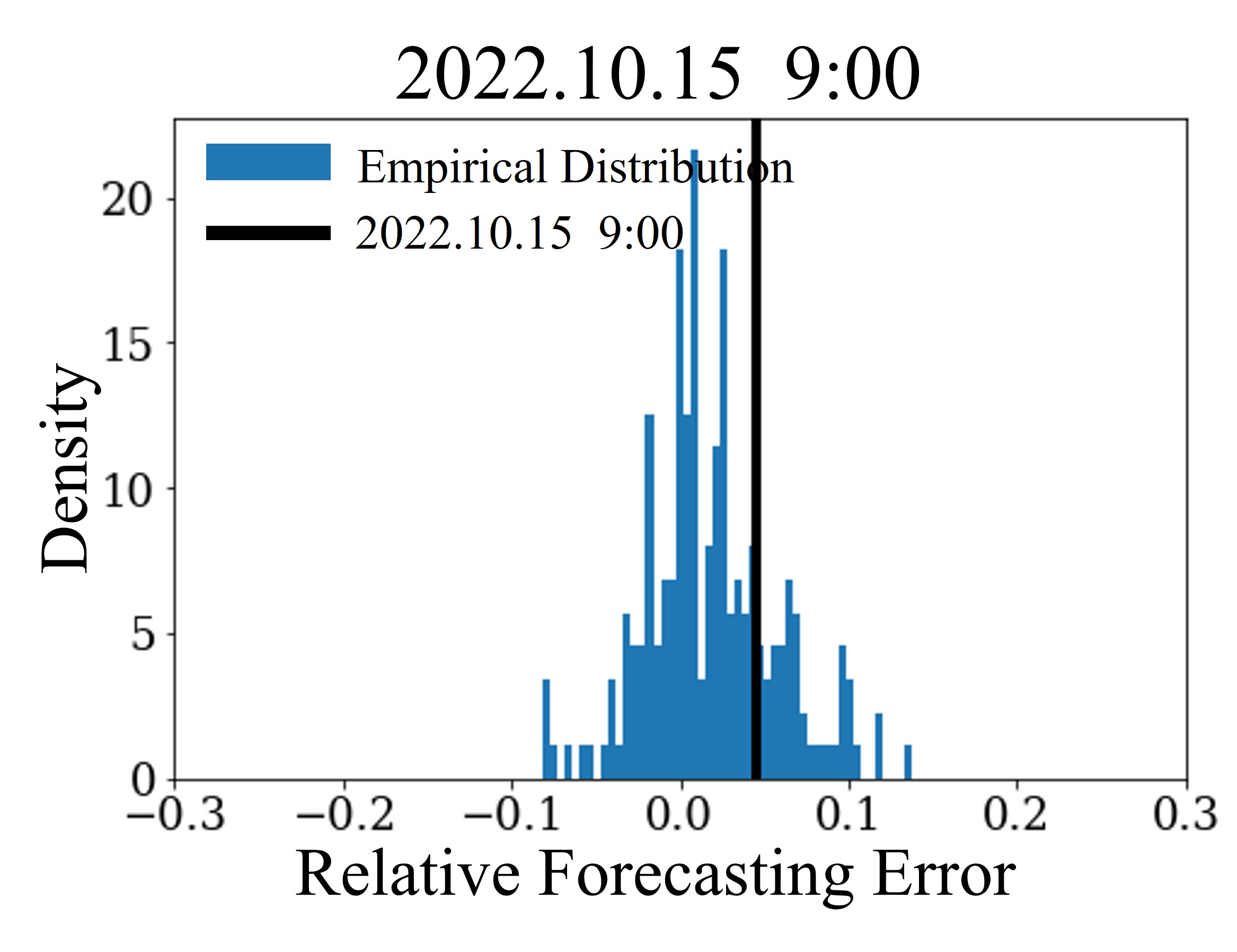}
  \caption{{\color{black}The \textit{future} error values and the scenarios' empirical distribution (proposed method)}}
  \label{new_distribution}
\end{figure}

\indent In this section, we mainly focus on the wind forecasting error affected by the \textit{deterministic} wind forecast value, the wind power ramp rate, temperature, and relative humidity. To illustrate the wind forecasting error under different wind forecast values, the empirical {\color{black}forecasting error} data is extracted from June to August 2022 in a similar geometric region in Texas. Fig.\ref{density_power}. shows the wind forecasting error density function under different normalized wind forecasting levels. It is clear that the empirical distribution under high wind output level ([0.8, 1], blue) is right-skewed to the lower wind output level ([0, 0.8], green), which means directly using the whole past scenario may make the decision aggressive under high wind weather. Similar results are also found in \cite{zhang2013modeling,mauch2013day}.
\indent Similarly, under different wind power ramp rates, the forecasting error has different patterns. The ramping rate of wind power is an aggregated parameter reflecting the stationary level of the environment. Typically, the high absolute wind power ramp rate means the weather data used for wind power forecasting is less precise than the stable environment, which results in different forecasting error distributions. After normalizing the ramping rate to [-1, 1] interval, the distribution of forecasting error conditions on a high ramp-up rate ([0.6, 1]) from June to August 2022 compared with its marginal distribution is shown in Fig.\ref{density_rate}.\\
\indent Unlike the wind power forecast value and ramp rate, the temperature and relative humidity have less influence on the wind power forecasting error. The distribution of wind power forecasting errors in southern Texas during hot and cold, dry and wet days of 2022 are compared in Fig.\ref{density_others}. It is obvious that the forecasting error under different temperatures or humidity shares a very similar distribution.
\indent To quantify the correlation level between the above four parameters and the wind power forecasting error, we calculate the Pearson correlation coefficient between each parameter and the forecasting error in southern Texas based on the past half year's data, i.e. February to July 2022. As shown in TABLE \ref{coefficient}, the wind power forecasting value and changing rate have more correlation with the forecasting error than the temperature and relative humility.
\begin{table}[H]
 \caption{The correlation coefficient between wind power forecasting error and some environmental parameters}
\centering
\begin{tabular}{cc}
\toprule
\textbf{Correlation between forecasting error and} & \textbf{Coefficient} \\ \midrule
Wind Power Forecasting Value                       & 0.50                 \\
(Absolute) Wind Power Changing Rate                & -0.17                \\
Temperature                                        & -0.07                \\
Relative Humility                                  & -0.03   \\ \bottomrule            
\label{coefficient}
\end{tabular}
\end{table}
\indent {\color{black}\textit{Remark}: Due to confidentiality concerns, the chosen four parameters are publicly available information accessible online. Decision-makers also have the option to incorporate locally measured data, including wind speed and direction, air pressure, and freezing level, to enhance the resolution of their analysis. By combining these diverse meteorological parameters, it becomes possible to delineate specific scenarios for various extreme weather events based on their respective thresholds.}\\
\indent After the above calculation, the scenarios under similar environments can be generated by Algorithm \ref{alg:samp} or \ref{alg:incre}, where the scenarios are selected from a smaller parameter space rather than the conventional method whose scenarios are randomly extracted from the whole past data. In other words, the distribution of empirical scenarios is adaptive changing with the environment, while it is fixed in the conventional method.\\
\begin{figure}[H]
\centering
  \includegraphics[scale=0.35]{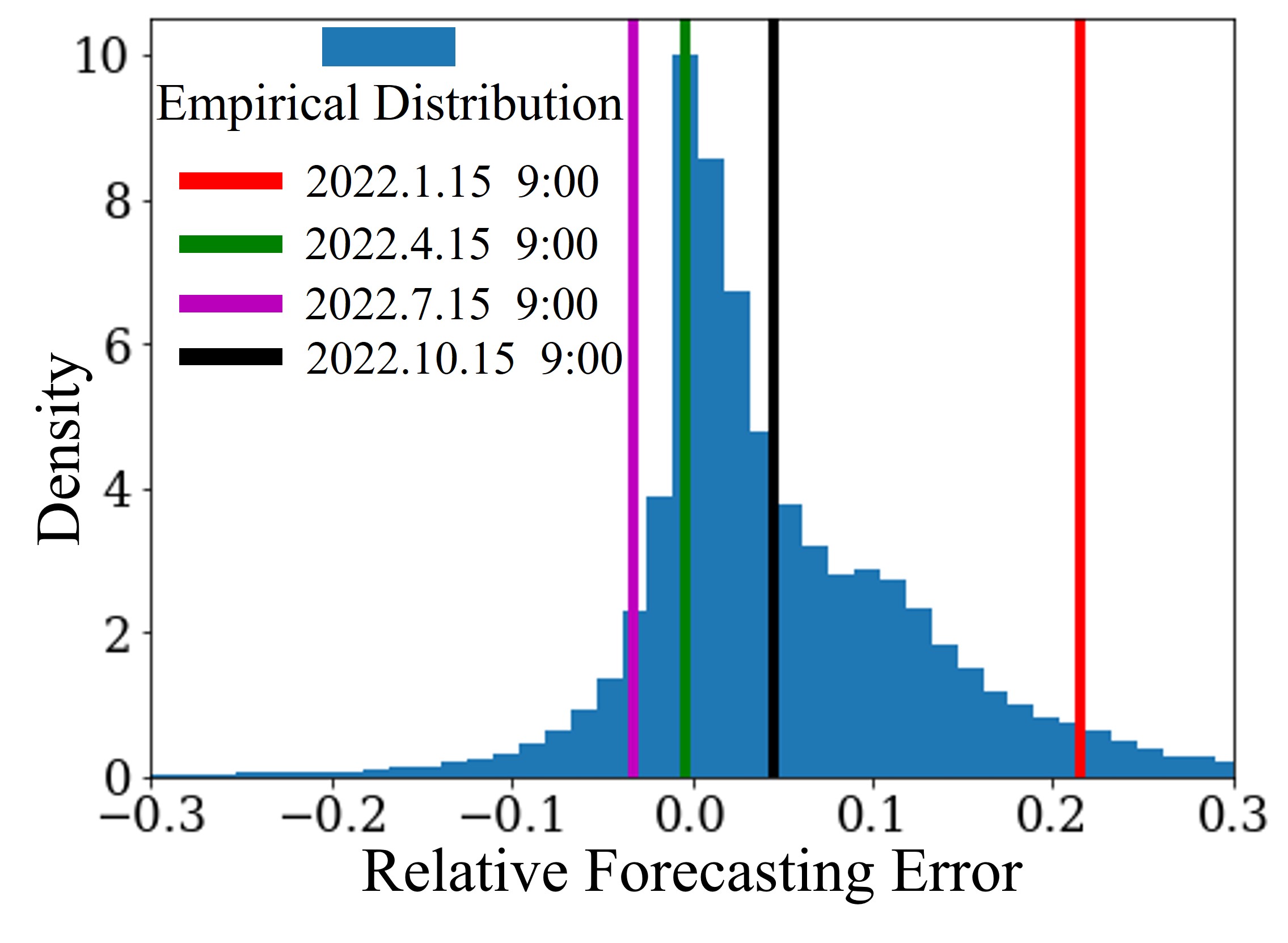}
  \caption{{\color{black}The \textit{future} error values and the scenarios' empirical distribution (conventional method)}}
  \label{old_distribution}
\end{figure}
\indent Supposing 200 scenarios are needed to meet the risk requirement, Fig.\ref{old_distribution} shows some real \textit{future} error values and the scenarios' empirical distribution in the whole last year, while the proposed method gives the adaptive empirical distribution for each decision-making time in Fig.\ref{new_distribution}.\\ 
\indent The advantages of sampling from similar parameter space are obvious from the above simulation. First, the distribution of scenarios will be compressed into narrower intervals giving a higher-resolution description of the future uncertainty. Second, the \textit{future} error values are also covered by this empirical distribution, which somewhat validates our method's correctness.
\subsection{24-bus System}
\begin{figure}[H]
\centering
  \includegraphics[scale=0.5]{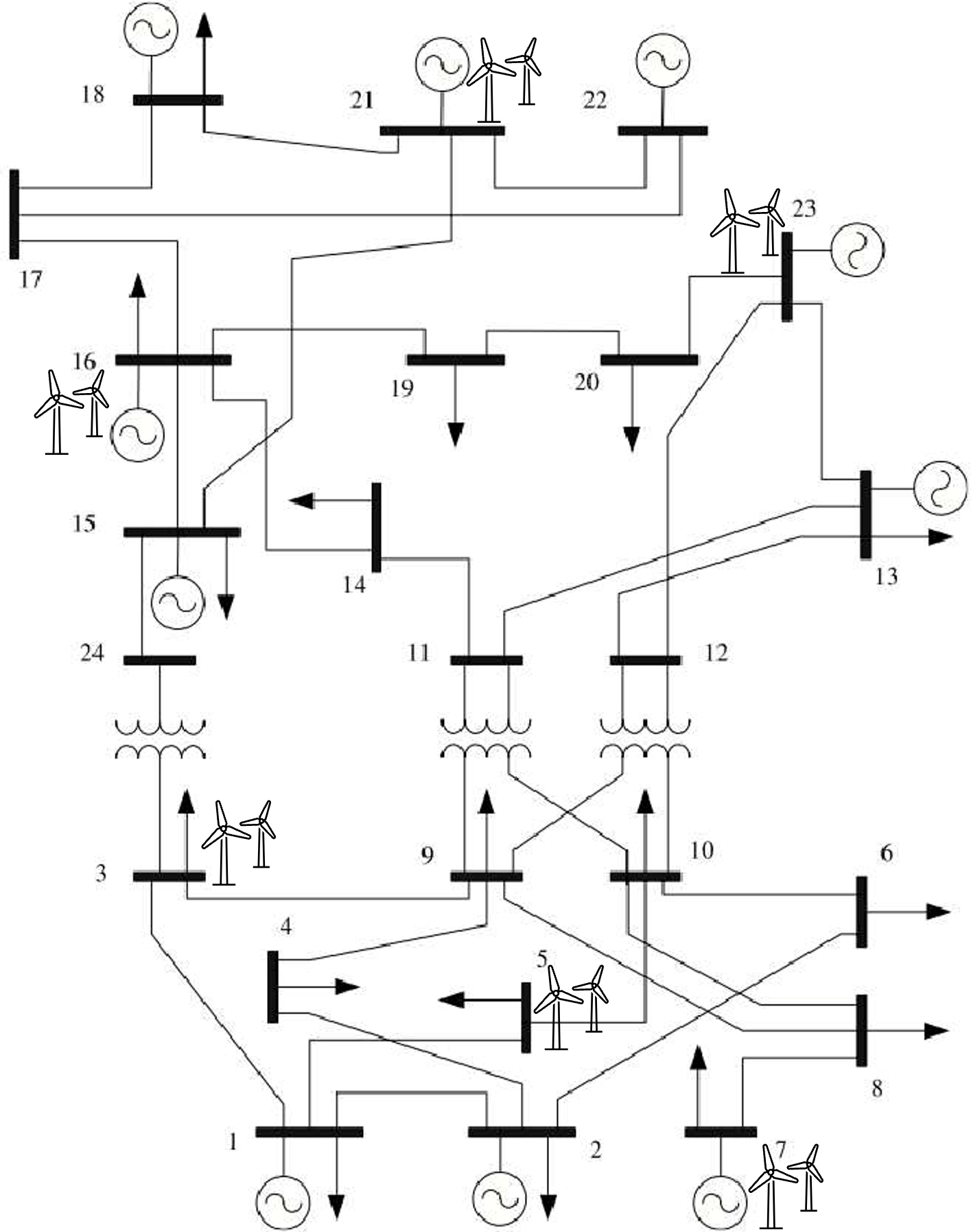}
  \caption{Modified 24-bus power system integrated with wind}
  \label{24bus}
\end{figure}

\begin{figure}[H]
\centering
  \includegraphics[scale=0.4]{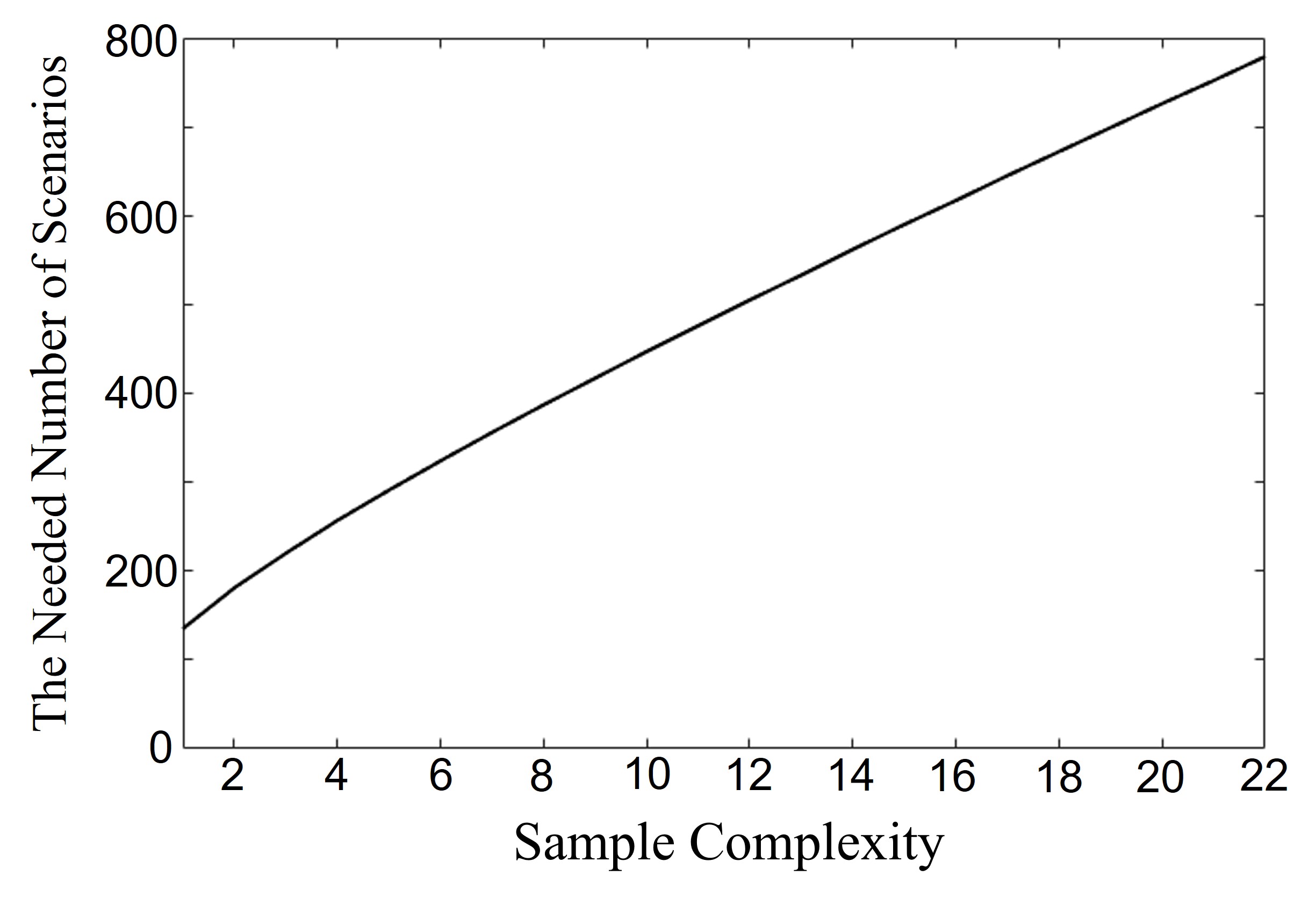}
  \caption{{\color{black}The relationship between the needed number of scenarios and sample complexity}}
  \label{complexity}
\end{figure}

\begin{figure}[H]
\centering
  \includegraphics[scale=0.35]{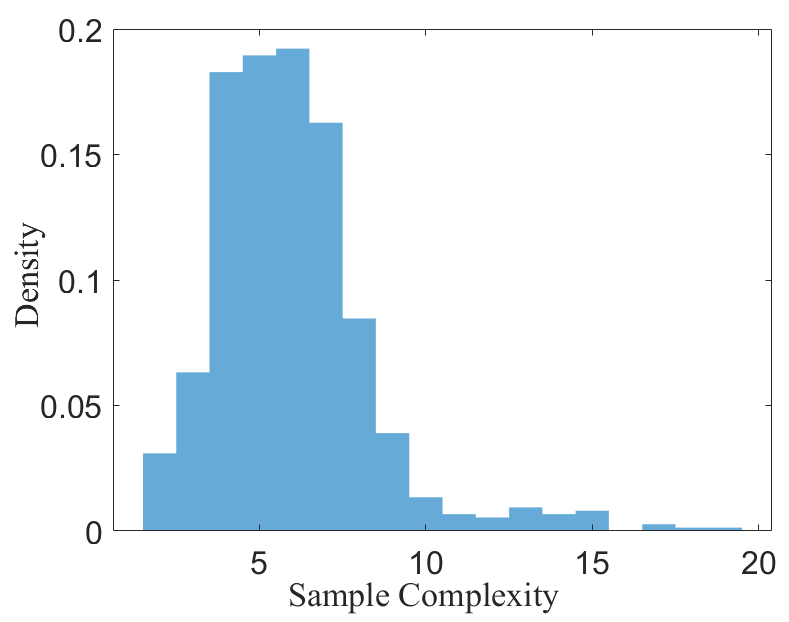}
  \caption{{\color{black}The density function of the sample complexity over the 744 dispatch intervals}}
  \label{complexitydensity}
\end{figure}

\begin{figure}[H]
\centering
  \includegraphics[scale=0.34]{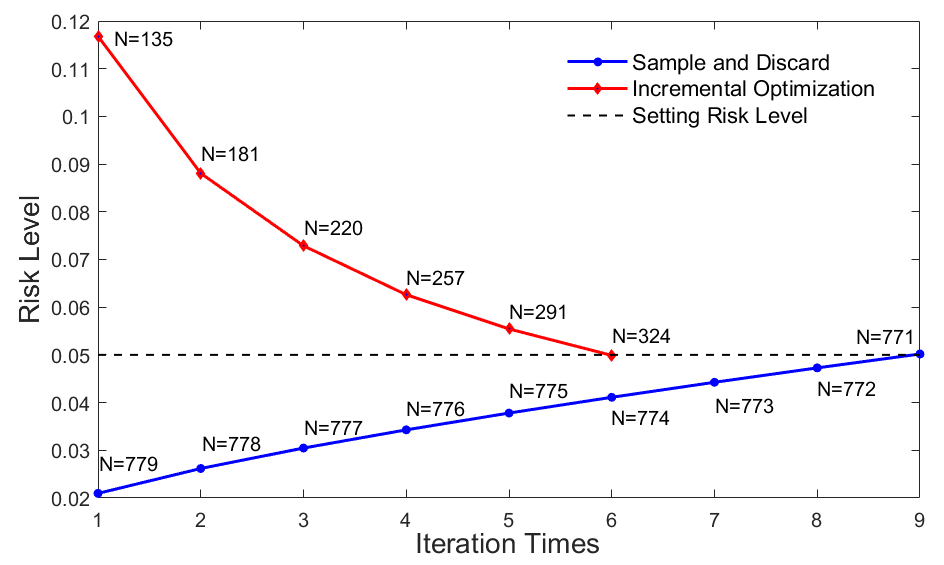}
  \caption{The risk tuning details of two different methods}
  \label{plottuning}
\end{figure}

\begin{figure}[H]
\centering
  \includegraphics[scale=0.3]{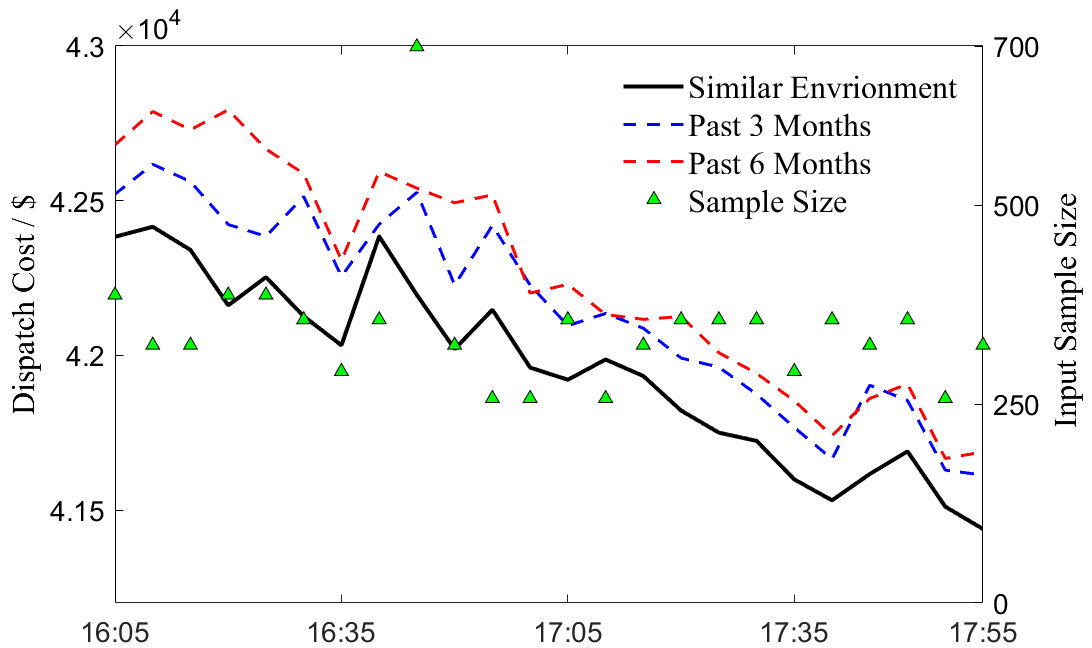}
  \caption{{\color{black}The dispatch cost by sampling from different parameter space and input sample size between 16:00 and 18:00 on August 2nd, 2022 (24-bus system, high wind day)}}
  \label{82}
\end{figure}

\begin{figure}[H]
\centering
  \includegraphics[scale=0.3]{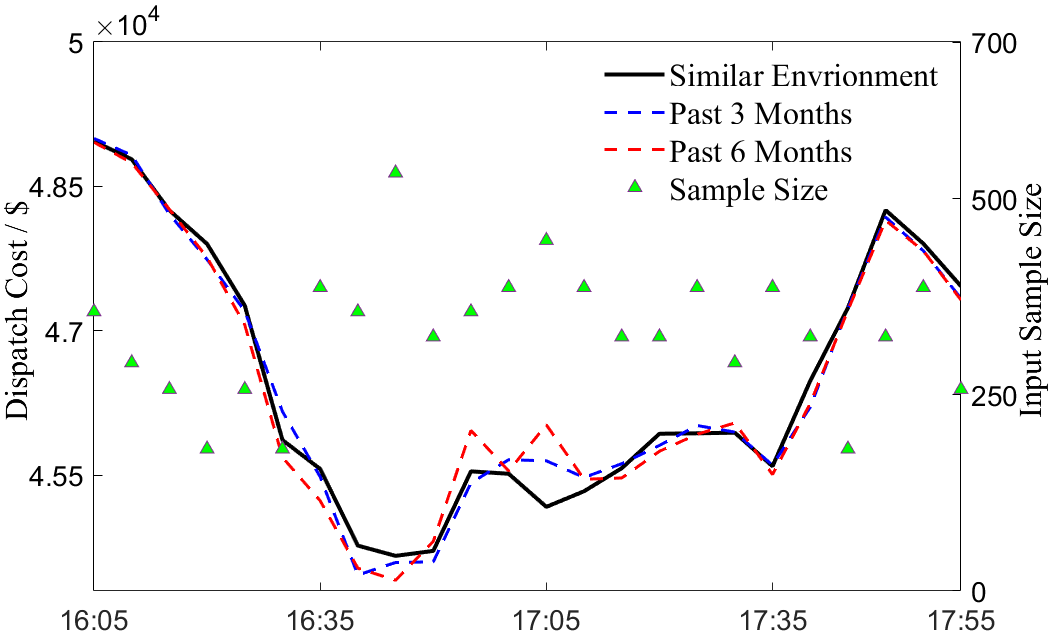}
  \caption{{\color{black}The dispatch cost by sampling from different parameter space and input sample size between 16:00 and 18:00 on August 30th, 2022 (24-bus system, low wind day)}}
  \label{830}
\end{figure}

\indent The 24-bus System is modified from the IEEE Reliability Test System (RTS-24) \cite{grigg1999ieee} with {\color{black}an additional} six wind farms to mimic the high renewable energy penetration as illustrated in Fig.\ref{24bus}. The detailed information, including generator parameters, reactance and capacity of transmission lines, and the load profile can be found in \cite{ordoudis2016updated}. The forecasting profile of wind generators located at bus 3, 5 and 7 are directly scaled from the west Texas region wind forecasting results, while the wind farms at bus 16, 21 and 23 are from the south Texas region. Each wind {\color{black}farm is} assumed with 400MW capacity and {\color{black}a low} marginal price (3 \$/MWh).\\
\indent As suggested in \cite{zugno2015robust}, the capacity on the transmission lines connecting the node pairs (15, 21), (14, 16), and (13, 23) is reduced to 400MW, 250MW, and 250MW, respectively. This is done to introduce bottlenecks or congestion in this high wind penetration system, which shares a similar situation when comes to the real power grids \cite{millstein2021solar}.\\
\indent The number of decision variables $n$ in the 24-bus system is 22 after eliminating the equation constraints. In this case, a rough approximation of Helly’s dimension $h$ by $n$ will increase the needed number of scenarios to meet the risk threshold, which results in great conservatism of the final decision. To illustrate this, Fig.\ref{complexity} shows the relationship between the needed number of scenarios and sample complexity under the risk $\epsilon=0.05$ and confidence parameter $\beta=10^{-3}$. It can be seen that when the sample complexity is 22, the needed number of scenarios is 779, which is much larger than the 324 needed scenarios under 6 sample complexity.
\indent Instead of using the a-prior approach, we apply Algorithm \ref{alg:incre} incrementally tuning risk for finding the exact sample complexity and the needed number of scenarios. Because the wind output and demand are changing, the sample complexity needs to be updated for each dispatch interval. Fig.\ref{complexitydensity} shows the empirical distribution of the sample complexity of the 744 dispatch intervals during August 16:00-18:00 we studied, where the average sample complexity (5.8) is much smaller than the number of decision variables (22). This evidence implies that in the previous \textit{sample and discard} risk-tuning method \cite{modarresi2018scenario}, where the sample complexity is first supposed to be the number of decision variables 22, the decision makers may need to discard some redundant scenarios to achieve the same performance as our incremental method.
\begin{table}[H]
\caption{Comparison of the risk tuning process}
\centering
\label{tuning}
\scalebox{0.85}{
\begin{tabular}{ccc}
\toprule
\textbf{Risk Tuning Method}        & \textbf{Sample and Discard \cite{campi2011sampling}} & \textbf{Incremental Optimization} \\ \midrule
Initial Input Sample Size & 779                         & 135                               \\
Intial Risk Level         & 0.021                       & 0.117                             \\
Final Input Sample Size   & 771                        & 324                               \\
Iteration Times           &  8
& 6            \\ \bottomrule                    
\end{tabular}}
\end{table}
\indent TABLE \ref{tuning} compares the initial and final input sample size of these two risk-tuning methods when sample complexity is 6. Compared with \textit{sample and discard}, incremental optimization guarantees the minimum input sample size with fewer iteration times, which is beneficial in the situation of limited high-accuracy data. Fig.\ref{plottuning} shows more details about the risk-tuning process of these two methods. The incremental optimization method's efficiency results in fewer iteration times and less sample size in each iteration.\\
\indent \textit{Remark}: The final input sample size in the sample and discard method is much larger than the incremental method because only support scenarios are discarded in \cite{campi2011sampling} but not randomly, and this risk bound is proved to be not tight in recent research \cite{romao2022exact}.
\indent {\color{black}So far, the simulation results have provided insights into the specifics of scenario generation and risk-tuning processes. Subsequent simulations will focus on examining the cost and actual violation outcomes following the scenario-based optimization model (\ref{saaa}).}\\
\indent After setting the same risk and confidence parameter as Fig.\ref{complexity}, we can input \textit{different} size of scenarios to meet the risk requirement based on the sample complexity of each dispatch interval, which is more efficient than the traditional sample and discarding method. Fig.\ref{82} illustrates the adaptive input sample size (green triangle) between 16:00 and 18:00 on August 2nd, 2022, where the dispatch cost of sampling from different parameter spaces is also compared \footnote{The dispatch cost is calculated after each dispatch with true wind output data, and the searching space of a similar environment is set to the past 3 months.}. During these two hours, the average wind power output is 23\% of the total wind generation capacity, which is a relatively high wind period during August rush hours. It is clear that sampling from a similar environment space results in lower costs in each dispatch interval than others.\\ 
\indent It needs to {\color{black}be clarified} that because of the stochastic property, sampling from a similar environment does not ensure a lower cost in each dispatch interval, especially on a low wind day. For example, the average wind power output is only 13\% of the total wind generation capacity between 16:00 and 18:00 on August 30th, which makes the economic benefits of sampling from a similar environment not comparable with the others (Fig.\ref{830}). \\
\indent However, sampling from a similar environment builds a more precise model of the uncertainty variables, which results in {\color{black}a} more trustworthy risk guarantee of the final solution. TABLE.\ref{result_24bus} compares the testing risk and the average cost of each 5-min dispatch interval by sampling from different parameter spaces during the rush hours 16:00 to 18:00 in August, where sampling from similar environment parameter space results in both lower testing violations and less dispatch cost. Meanwhile, when randomly sampling from the past half year, the testing violation (0.054) even exceeds the setting risk (0.05), which invalidates the risk-guarantee property of the scenario approach.
\begin{table*}[htbp]
\caption{The violation and average dispatch cost of sampling from different parameter spaces (24-bus system)}
\centering
\label{result_24bus}
\begin{tabular}{ccccc}
\toprule
&\textbf{Sampling Space}       & \textbf{Past Half Year} & \textbf{Past Three Months} & \textbf{Similar Environment} \\\midrule
\multicolumn{1}{c}{\multirow{5}{*}{\textbf{Whole August}}} & {\color{black}\textbf{Setting Violation}}            & {\color{black}0.05}                    & {\color{black}0.05}                     & {\color{black}0.05}  \\
& \textbf{Actual Violation}            & \textbf{0.054}                    & 0.040                     & 0.036  \\
\multicolumn{1}{c}{} &\textbf{Average Cost($10^4\$$)} & 4.7522                  & 4.7512                     & 4.7444  \\
\multicolumn{1}{c}{} &{\color{black}\textbf{Average Solving Time(s)}} & {\color{black}0.181}                  & {\color{black}0.166}                     & {\color{black}0.165}  \\
\multicolumn{1}{c}{} &{\color{black}\textbf{Average Sampling Time(s)}} & {\color{black}0.015}                  & {\color{black}0.014}                     & {\color{black}0.020}  \\\midrule       \textbf{2nd August (high wind)}   & \textbf{Average Cost($10^4\$$)}            & 4.2257                    & 4.2154                       & 4.1956        \\ 
\textbf{30th August (low wind)}                     & \textbf{Average Cost($10^4\$$)} & 4.6414                 & 4.6430                     & 4.6433 \\ \bottomrule
\end{tabular}
\end{table*}
\begin{table*}[h]
\caption{The violation and average dispatch cost of sampling from different parameter spaces (118-bus system)}
\centering
\label{result_118bus}
\begin{tabular}{ccccc}
\toprule
&\textbf{Sampling Space}       & \textbf{Past Half Year} & \textbf{Past Three Months} & \textbf{Similar Environment} \\\midrule
\multicolumn{1}{c}{\multirow{5}{*}{\textbf{Whole August}}} & {\color{black}\textbf{Setting Violation}}            & {\color{black}0.05}                    & {\color{black}0.05}                     & {\color{black}0.05}  \\
& \textbf{Actual Violation}            & 0.042                    & 0.035                     & 0.035  \\
\multicolumn{1}{c}{} &\textbf{Average Cost($10^5\$$)} & 1.7445                  & 1.7454                     & 1.7414  \\
\multicolumn{1}{c}{} &{\color{black}\textbf{Average Solving Time(s)}} & {\color{black}13.66}                  & {\color{black}13.41}                     & {\color{black}13.30}  \\
\multicolumn{1}{c}{} &{\color{black}\textbf{Average Sampling Time(s)}} & {\color{black}0.022}                  & {\color{black}0.021}                     & {\color{black}0.028}  \\\midrule       \textbf{2nd August (high wind)}   & \textbf{Average Cost($10^5\$$)}            & 1.5294                    & 1.5435                       & 1.5173        \\ 
\textbf{30th August (low wind)}                     & \textbf{Average Cost($10^5\$$)} & 1.7251                  & 1.7255                     & 1.7271 \\ \bottomrule                
\end{tabular}
\end{table*}

\subsection{118-bus System}
\begin{figure}[H]
\centering
  \includegraphics[scale=0.3]{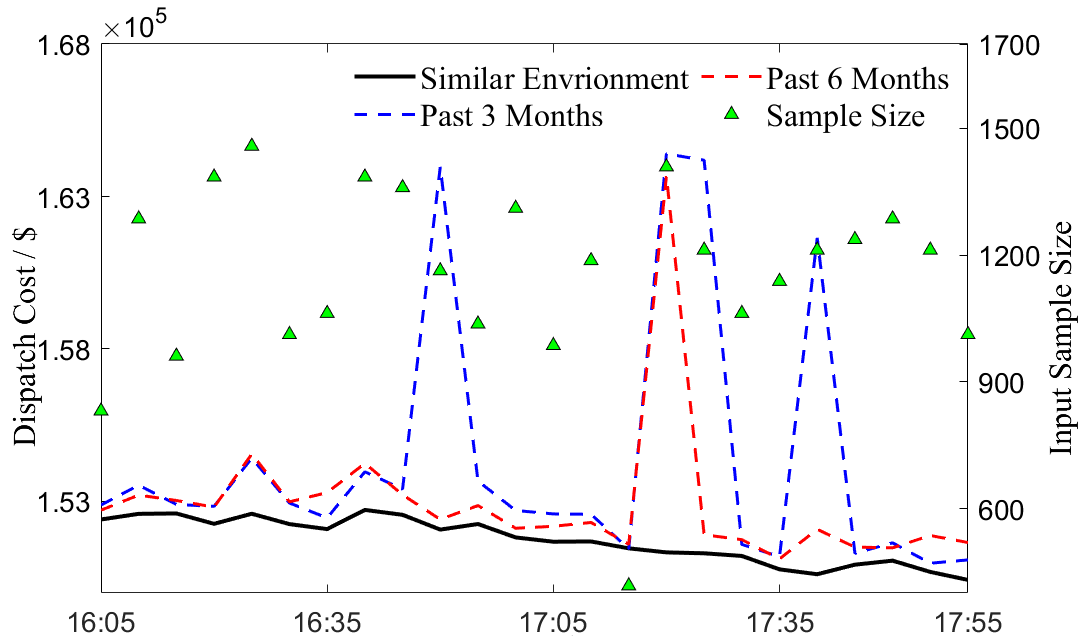}
  \caption{{\color{black}The dispatch cost by sampling from different parameter space and input sample size between 16:00 and 18:00 on August 2nd, 2022 (118-bus system, high wind day)}}
  \label{118-82}
\end{figure}

\indent We apply our method to a larger system in this section, based on the test case \textit{c118swf.m} in MATPOWER \cite{murillo2013secure}. This system includes 118 nodes, 210 lines, and 52 generators, 11 of which are modeled as wind farms. To address the influence of wind uncertainty on the economic dispatch, we replace the 4 storage units with the same capacity wind farms. Meanwhile, the transmission line capacities are set to be 60\% of the original value to introduce more congestion in the system.\\
\indent Similar to the simulation in the 24-bus system, the wind power forecasting and real value are directly modeled by the data from ERCOT with 7 wind generators from south Texas and 8 wind generators from west Texas. Fig.\ref{118-82} shows the adaptive input sample size and the dispatch cost of sampling from different parameter spaces between 16:00 and 18:00 on August 2nd, 2022. In this period, the average wind power output accounts for 24\% of the total wind generation capacity, a relatively high wind output level during August peak hours.\\
\indent Benefiting from the smaller uncertainty space (see Fig.\ref{old_distribution}, Fig.\ref{new_distribution}), sampling from a similar environment helps the decision maker avoid some odd scenarios when applying scenario approach in the real world. These odd scenarios typically make the final decision more conservative, which can be illustrated by the cost spikes in Fig.\ref{118-82} when we directly sample scenarios from the past.
\indent TABLE.\ref{result_118bus} lists the testing risk and the average cost of 5-min dispatch during August peak hours (16:00-18:00), where the result has a similar pattern with the 24-bus system. The 5-min short-term wind power forecasting error is relatively lower than other long-term forecasting errors, but a more efficient scenarios generation model will make the chance-constrained dispatch solution more trustworthy.\\
\indent Although there's no theory guarantee for a lower cost solution by sampling from a similar environment, the simulation shows the economic advantage of the proposed efficient scenarios generation model, especially during high wind situations. To quantify the relationship between the accuracy of the scenario generation model and economic benefits will be one possible direction of future research.\\
\indent {\color{black}Regardless of the selection of sample space, the number of needed scenarios is the same, which results in a similar problem-solving time of three different sample spaces in both TABLE.\ref{result_24bus} and TABLE.\ref{result_118bus}. Meanwhile, computing \textit{indicator vector} (\ref{correlation3}) of similar environments increases the sampling time, but it's far less than the solving time by one or more orders of magnitude.}
\section{Conclusion}\label{sec:conc}
\indent This paper studies the main two barriers to applying the scenario approach to the economic dispatch with high penetration of renewable resources, i.e. lack of accurate scenario generation models and inefficient risk tuning process. Leveraging correlation analysis, we generate scenarios via an environment filter with empirical distribution closer to the true probability measure. After embedding this scenario generation model with the incremental scenario optimization algorithm, we propose an efficient risk-tuning scheme, which can solve the optimal solution meeting risk requirement with minimum data size and provide other higher-risk solutions to system operators in meanwhile. Case studies based on real-world wind data and modified IEEE benchmark systems show the effectiveness and advantages of our methods.\\
\indent Directly generating scenarios from past experience may be a naive approach, but it works well in practice when the needed data size is much smaller than the size of the past data pool. Future work includes (1) comparing the results of using other scenario generation methods, such as the generative model in the field of machine learning; (2) extending the economic dispatch problem to a multi-stage framework, i.e. look-ahead economic dispatch, {\color{black}where varying risk tolerance levels are allowed across different time horizons to enhance the adaptability}; and (3) applying the proposed scheme to unit commitment and other non-convex decision-making processes in electric power systems, {\color{black} which requires efficient handling of discrete variables under incremental scenarios to ensure timely solutions}.

\bibliographystyle{IEEEtran}
\bibliography{ref.bib}

\begin{IEEEbiography}[{\includegraphics[width=1in,height=1.25in,clip,keepaspectratio]{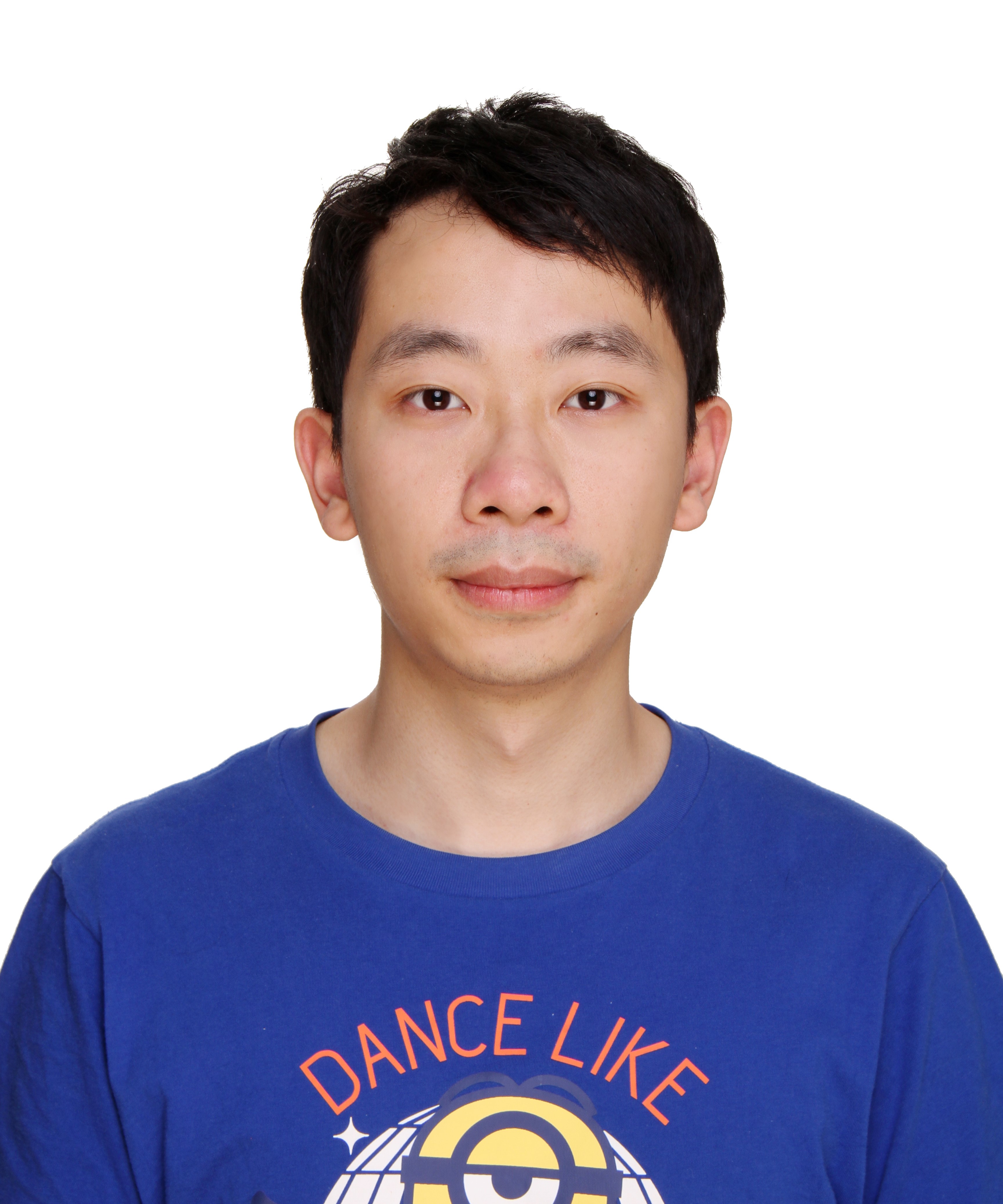}}]{Qian Zhang}
(Student Member, IEEE) received his B.E. and M.S. degrees in Electrical Engineering from Zhejiang University, Hangzhou, China, in 2019 and 2022 respectively. He is currently working toward the Ph.D. degree at Texas A\&M University, College Station, TX, USA. His research interests include machine learning, optimization in the electricity market, and power system stability and control.  
\end{IEEEbiography}

\begin{IEEEbiography}[{\includegraphics[width=1in,height=1.25in,clip,keepaspectratio]{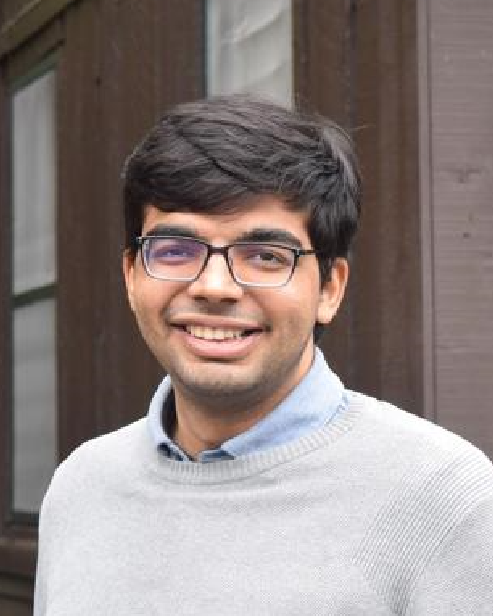}}]{Apurv Shukla}
(Member, IEEE) received his bachelor’s degree in Manufacturing Science and Engineering from IIT Kharagpur in 2016 and his Ph.D. degree in Operations Research from Columbia University in 2022. He is a Postdoctoral Associate at Texas A\&M University. His research interest lies in systems and control.
\end{IEEEbiography}

\begin{IEEEbiography}[{\includegraphics[width=1in,height=1.25in,clip,keepaspectratio]{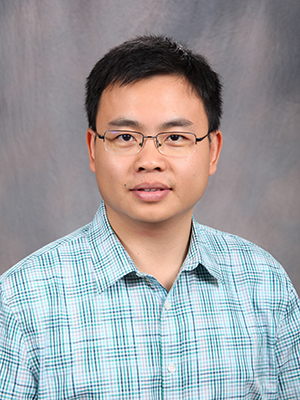}}]{Le Xie}
(Fellow, IEEE) received the B.E. degree in electrical engineering from Tsinghua University, Beijing, China, in 2004, the M.S. degree in engineering sciences from Harvard University, Cambridge, MA, USA, in 2005, and the Ph.D. degree from the Department of Electrical and Computer Engineering, Carnegie Mellon University, Pittsburgh, PA, USA, in 2009. He is currently a Professor with the Department of Electrical and Computer Engineering, Texas A\&M University, College Station, TX, USA. His research interests include modeling and control of large-scale complex systems, smart grid applications with renewable energy resources, and electricity markets.
\end{IEEEbiography}

\end{document}